\documentclass{article}
\usepackage{amsmath}
\usepackage{amssymb}
\usepackage{amsthm}
\usepackage{mathtools}
\usepackage{tikz-cd}
\usepackage{enumitem}
\usepackage{authblk}

\usepackage[backend=biber,maxbibnames=9,giveninits=true,isbn=false]{biblatex}
\AtEveryBibitem{\clearfield{number}}

\newtheorem{theorem}{Theorem}
\newtheorem{lemma}[theorem]{Lemma}
\newtheorem{proposition}[theorem]{Proposition}
\newtheorem{corollary}[theorem]{Corollary}
\newtheorem{definition}[theorem]{Definition}
\newtheorem{example}[theorem]{Example}
\theoremstyle{remark}
\newtheorem{remark}[theorem]{Remark}

\title{A Characterization of Circle Graphs\\ in Terms of Total Unimodularity}
\author[1]{Robert Brijder}
\author[2]{Lorenzo Traldi}
\affil[1]{Hasselt University, Belgium}
\affil[2]{Lafayette College, Easton, PA, USA}
\date{}

\DeclareMathOperator{\spn}{span}
\DeclareMathOperator{\Tch}{Tch}

\newcommand{\dTch}[2]{\Tch_{#2}(#1)}
\newcommand{\edgetrans}[2]{\mathrm{eti}_{#1,#2}}
\newcommand{\eCnt}[2]{\sigma(#1,#2)}

\begin{document}

\maketitle

\begin{abstract}
A graph $G$ has an associated multimatroid $\mathcal{Z}_3(G)$, which is equivalent to the isotropic system of $G$ studied by Bouchet. In previous work it was shown that $G$ is a circle graph if and only if for every field $\mathbb F$, the rank function of $\mathcal{Z}_3(G)$ can be extended to the rank function of an $\mathbb F$-representable matroid. In the present paper we strengthen this result using a multimatroid analogue of total unimodularity. As a consequence we obtain a characterization of matroid planarity in terms of this total-unimodularity analogue.
\end{abstract}

\section{Introduction}

The outline of the theory of circle graphs and local complementation was set forth by Andr\'e Bouchet in a series of papers published over several decades. Much of his work involved two kinds of combinatorial structures, delta-matroids \cite{bouchet1987, mp/Bouchet87} and isotropic systems \cite{Bouchet/87/ejc/isotropicsys}. In the late 1990s, Bouchet unified these two structures by introducing a common generalization, called a \emph{multimatroid} \cite{DBLP:journals/siamdm/Bouchet97}. 

To state the definition of a multimatroid we need some terminology regarding partitions. If $\Omega$ is a partition of a set $U$ then the elements of $\Omega$ are called \emph{skew classes}. A \emph{transversal} of $\Omega$ is a subset of $U$ that contains precisely one element of every skew class, and a \emph{subtransversal} of $\Omega$ is a subset of a transversal. The sets of subtransversals and transversals of $\Omega$ are denoted $\mathcal{S}(\Omega)$ and $\mathcal{T}(\Omega)$, respectively. We use $2^S$ to denote the power set of a set $S$. In this paper we assume knowledge of some basic matroid-theoretic notions, see, e.g., \cite{Oxley/MatroidBook-2nd,Oxley/MatroidTutor/2003} for an introduction.

\begin{definition} \label{def:semi-multimatroid}
A \emph{multimatroid} $Z$ (described by its rank function) is a triple $(U,\Omega,r)$, where $\Omega$ is a partition of a finite set $U$ and $r:\mathcal{S}(\Omega) \to \mathbb N$ is a function such that for each $S \in \mathcal{S}(\Omega)$
\begin{itemize}
\item $(S,r|_{2^S})$ is a matroid (described by its rank function), where $r|_{2^S}$ denotes the restriction of $r$ by $2^S$, and
\item if $x$ and $y$ are distinct elements of a skew class of $\Omega$ disjoint from $S$, then $\max \{ r(S\cup \{x\}), r(S\cup \{y\}) \} > r(S)$.
\end{itemize}
\end{definition}
If $T \in \mathcal{T}(\Omega)$, then the matroid $(T,r|_{2^T})$, denoted by $Z[T]$, is called the \emph{transverse matroid} of $Z$ corresponding to $T$. Also, if each skew class has at least two elements, then $Z-T := (U-T,\Omega',r|_{2^{U-T}})$ with $\Omega' = \{ \omega - T \mid \omega \in \Omega \}$ is a multimatroid. 

A multimatroid in which every skew class has exactly $k$ elements is called a \emph{$k$-matroid}.

\begin{definition}
A multimatroid $Z=(U,\Omega,r)$ is \emph{sheltered} by a matroid $M$ if $M$ is a matroid on $U$ whose rank function restricts to $r$. If the rank of $M$ is the maximum value of $r(S)$ with $S \in \mathcal{S}(\Omega)$, then $M$ is a \emph{strict} sheltering matroid for $Z$.
\end{definition}

In his fourth paper on multimatroids \cite{Bouchet/MMIV}, Bouchet introduced a rather complicated notion of representability inspired by both his notion of representability for delta-matroids \cite{bouchet1987} and Tutte's notion of matroid representability using chain groups. We do not know of any other research on multimatroids that has been done using Bouchet's notion of representability.

In recent work a different notion of multimatroid representability is used, which seems more natural: a (strict) $\mathbb F$-\emph{representation} of a multimatroid $Z$ is an $\mathbb F$-representation of a (strict) sheltering matroid for $Z$. We say that $Z$ is (strictly) \emph{representable} over $\mathbb F$ if it has a (strict) $\mathbb F$-representation.

We use the following notation for matrices. If $X$ and $Y$ are finite sets then an $X \times Y$ matrix has rows and columns that are not ordered, but are indexed by $X$ and $Y$ (respectively). Suppose $G$ is a \emph{looped simple graph}, i.e., a graph which may have loops but has no more than one loop at any vertex, and no more than one edge connecting any two vertices. The \emph{adjacency matrix} $A(G)$ of $G$ is a $V(G) \times V(G)$-matrix over $GF(2)$, where, for $u,v \in V(G)$, the entry of $A(G)$ indexed by $(u,v)$ is $1$ if and only if there is an edge between $u$ and $v$. In particular, loops are represented by nonzero diagonal entries in $A(G)$. Recall that an $X \times Y$ matrix $A$ \emph{represents} a matroid $M$ with ground set $Y$, where, for all $Y' \subseteq Y$, the rank of $Y'$ in $M$ is equal to the rank of $A$ restricted to the columns of $Y'$. A matroid is called \emph{binary} if it is represented by a matrix over $GF(2)$.

\begin{definition}
If $G$ is a looped simple graph, then the \emph{isotropic matroid} $M[\mathrm{IAS}(G)]$ of $G$ is the binary matroid represented by the $GF(2)$-matrix
\[
\mathrm{IAS}(G)=
\begin{pmatrix}
I & A(G) & I+A(G)
\end{pmatrix}
\text{,}
\]
where $I$ is the $V(G) \times V(G)$ identity matrix.
\end{definition}

Each vertex $v \in V(G)$ corresponds to a 3-element subset of the ground set $U$ of $M[\mathrm{IAS}(G)]$, called the \emph{vertex triple} of $v$, consisting of the column indices corresponding to the $v$-columns of $I$, $A(G)$, and $I+A(G)$ in $\mathrm{IAS}(G)$. Note that the vertex triples partition $U$. It turns out that $M[\mathrm{IAS}(G)]$ shelters a $3$-matroid $(U, \Omega, r)$, denoted by $\mathcal{Z}_3(G)$, where $\Omega$ is the set of vertex triples of vertices of $G$. 
Notice that $\mathrm{IAS}(G)$ provides a strict $GF(2)$-representation of $\mathcal{Z}_3(G)$, so $\mathcal{Z}_3(G)$ is strictly representable over $GF(2)$.

Recall that a matroid $M$ is \emph{regular} if it satisfies any of these equivalent conditions. (See, e.g., \cite{Oxley/MatroidBook-2nd}.)
\begin{enumerate}[label=(\alph*)]
\item \label{lbl:reg:GF2_diff} $M$ is representable over $GF(2)$ and some field of characteristic $\neq 2$.

\item \label{lbl:reg:all_fields} $M$ is representable over all fields.

\item \label{lbl:reg:TU} $M$ is represented over $\mathbb R$ by a matrix of integers $U$ which is \emph{totally unimodular}, i.e., every square submatrix of $U$ has determinant in the set $\{-1,0,1\}$.
\end{enumerate}
The smallest non-regular binary matroids are the Fano matroid $F_7$ and its dual $F_7^*$, represented by the $GF(2)$-matrices $(I_3 \enskip A)$ and $(I_4 \enskip A^T)$ respectively, where $I_3$ and $I_4$ are the identity matrices of dimensions $3$ and $4$ and
\[
    A =
    \begin{pmatrix}
    0 & 1 & 1 & 1 \\
    1 & 0 & 1 & 1 \\
    1 & 1 & 0 & 1
    \end{pmatrix}.
\]

At first glance the theory of regular matroids does not seem to be relevant to isotropic matroids. For instance, the above representation of $F_7$ is a submatrix of $\mathrm{IAS}(K_3)$, and so the isotropic matroid of $K_3$ is not regular. The next result from \cite{BT/IsotropicMatIV} shows however that the multimatroids associated with circle graphs have some special properties reminiscent of regular matroids. We recall the definition of a circle graph in Section~\ref{sec:cycle_matrices}. Recall that $\mathcal{Z}_3(G)$ is, by definition, representable over $GF(2)$.

\begin{theorem}[\cite{BT/IsotropicMatIV}]
\label{prequel}
These properties of a simple graph $G$ are equivalent.
\begin{enumerate}
\item \label{item:isom_circle_graph} $G$ is a circle graph.

\item \label{item:transv_unim} $\mathcal{Z}_3(G)$ has a strict representation $A$ over $\mathbb R$ containing only integer entries that is ``transversely unimodular''. That is, for every transversal $T$ of the set of skew classes of $\mathcal{Z}_3(G)$, the determinant of the square submatrix obtained from $A$ by retaining only the columns of $T$ is in $\{-1,0,1\}$.

\item \label{item:repr_GF2_diff} $\mathcal{Z}_3(G)$ is representable over some field of characteristic different from $2$.

\item \label{item:min_T_reg} For every transversal $T$ of the set of skew classes of $\mathcal{Z}_3(G)$, the $2$-matroid $\mathcal{Z}_3(G)-T$ is representable over some field of characteristic different from $2$.
\end{enumerate}
\end{theorem}

Properties~\ref{item:repr_GF2_diff} and \ref{item:min_T_reg} remain equivalent to the others if the phrase ``over some field of characteristic different from $2$'' is replaced with ``over all fields''. These equivalences are strongly reminiscent of the equivalent descriptions \ref{lbl:reg:GF2_diff} and \ref{lbl:reg:all_fields} of regular matroids mentioned above. On the other hand, property~\ref{item:transv_unim} of Theorem~\ref{prequel} seems weaker than the analogous property \ref{lbl:reg:TU} of regular matroids, as the unimodularity property of property~\ref{item:transv_unim} applies only to submatrices corresponding to transversals, not arbitrary subtransversals. In fact, square submatrices corresponding to subtransversals in the representation matrices considered in \cite{BT/IsotropicMatIV} can have various determinants; for instance, some entries of these matrices are equal to $2$. It is important to realize that property~\ref{item:repr_GF2_diff} of Theorem~\ref{prequel} does not require $\mathcal{Z}_3(G)$ to have a single sheltering matroid that is representable both over $GF(2)$ and over some field of characteristic $\neq 2$; there may be different sheltering matroids representable over different fields (see Remark~\ref{remark:no_TU_repr} below).

It is also important to realize that property~\ref{item:min_T_reg} of Theorem~\ref{prequel} implies that for every transversal $T$, the transverse matroid $\mathcal{Z}_3(G)[T]$ is regular; but this property is strictly weaker than property~\ref{item:min_T_reg}. For instance, it is easy to see that even though the wheel graph $W_5$ is not a circle graph, the transverse matroids of $\mathcal{Z}_3(W_5)$ are all regular. Indeed, the smallest non-regular binary matroids $F_7$ and $F^*_7$ each have $7$ elements, while the transverse matroids of $\mathcal{Z}_3(W_5)$ have only $6$ elements.

The implications $2 \implies 3$ and $3 \implies 4$ of Theorem~\ref{prequel} are fairly obvious, and $4 \implies 1$ is a fairly direct consequence of Bouchet's well-known characterization of circle graphs by forbidden vertex-minors. The difficult part of the proof of Theorem~\ref{prequel} in \cite{BT/IsotropicMatIV} is a long and technical argument that verifies the implication $1 \implies 2$ using interlacement graphs with respect to Euler systems in $4$-regular graphs. 

In the present paper we strengthen the proof of the implication $1 \implies 2$ in Theorem~\ref{prequel}. Let us say that a representation of $\mathcal{Z}_3(G)$ over $\mathbb R$ is \emph{totally transversally unimodular} if every square submatrix of that representation such that the column indices form a subtransversal has determinant in $\{-1,0,1\}$. We show the following (using Theorem~\ref{prequel} for the if direction).

\begin{theorem}\label{thm:main_result}
Let $G$ be a simple graph. Then $G$ is a circle graph if and only if the $3$-matroid $\mathcal{Z}_3(G)$ has a totally transversally unimodular representation. If this is the case then there exists a totally transversally unimodular representation of $\mathcal{Z}_3(G)$ that is strict.
\end{theorem}

Since the requirement of total transversal unimodularity applies to all subtransversals, Theorem~\ref{thm:main_result} provides a property analogous to property~\ref{lbl:reg:TU} of regular matroids. This completes the analogy between regular matroids and circle graphs. The proof of Theorem~\ref{thm:main_result} also provides a new insight into the situation by highlighting a natural connection, through a new notion called the edge-transition incidence matrix, between the cycle space of a 4-regular graph $F$, the cycle spaces of touch-graphs of circuit partitions of $F$, and the transverse matroids of $\mathcal{Z}_3(G)$ of circle graphs $G$ associated to Euler systems of $F$. (We recall the definition of a touch-graph in Section~\ref{sec:touch_graphs}.)

A matroid is said to be \emph{planar} if it is isomorphic to the cycle matroid of a planar graph. Recall that the \emph{fundamental graph} of a matroid $M$ with respect to a basis $B$ of $M$ is the bipartite graph where for $x \in B$ and $y$ an element of $M$ not in $B$, $x$ is adjacent to $y$ if $x$ is in the fundamental circuit of $y$ with respect to $B$. In \cite[Theorem~50]{BT/IsotropicMatIV} it is shown that a binary matroid $M$ is planar if and only if $\mathcal{Z}_3(G)$, with $G$ a fundamental graph of $M$, is representable over $GF(2)$ and over some field of characteristic different from $2$. Since this property of $\mathcal{Z}_3(G)$ is equivalent to $\mathcal{Z}_3(G)$ having a totally transversally unimodular representation, we immediately obtain the following characterization of planarity.
\begin{corollary}
Let $M$ be a binary matroid. Then the following three conditions are equivalent:
\begin{itemize}
\item $M$ is planar,

\item the $3$-matroid $\mathcal{Z}_3(G)$ has a strict, totally transversally unimodular representation for some fundamental graph $G$ of $M$, and

\item the $3$-matroid $\mathcal{Z}_3(G)$ has a strict, totally transversally unimodular representation for every fundamental graph $G$ of $M$.
\end{itemize}
\end{corollary}

Finally, in Section~\ref{sec:cycle_basis_Euler} we efficiently obtain new proofs of some essential results of \cite{BT/IsotropicMatIV} using constructions introduced in this paper.

\section{Preliminaries}\label{sec:prelim}

The main purpose of this section is to fix definitions of some well-known graph-theoretical notions.

\subsection{Walks and circuits}

We consider graphs where loops and multiple edges are allowed. The notion of a half-edge will be important in this paper and so we explicitly define graphs using half-edges.

A \emph{graph} $G$ is a $4$-tuple $(V,H,E,\epsilon)$, where $V$ and $H$ are finite sets, $E$ is a partition of $H$ in (unordered) pairs, and $\epsilon: H \to V$ is a function. The elements of $V$, $H$, and $E$ are called \emph{vertices}, \emph{half-edges}, and \emph{edges} of $G$, respectively. We denote $V$, $H$, and $E$ by $V(G)$, $H(G)$, and $E(G)$, respectively. The number of connected components of $G$ is denoted by $c(G)$.

A \emph{directed graph} is defined analogously; the only difference is that $E$ is then a partition of $H$ in \emph{ordered} pairs. In that case, for $e = (h_1,h_2) \in E$, $h_1$ and $h_2$ are called the \emph{tail} and \emph{head} of $e$, respectively. We also say that $e$ is directed \emph{from} $\epsilon(h_1)$ \emph{toward} $\epsilon(h_2)$.

A vertex $v$ and half-edge $h$ are called \emph{incident} if $\epsilon(h) = v$. A \emph{single transition} is an unordered pair $\{h_1,h_2\}$ of half-edges incident to a common vertex. A \emph{directed} single transition is an ordered pair $(h_1,h_2)$ of half-edges incident to a common vertex; we say that $h_1$ is directed \emph{toward} the vertex, and $h_2$ is directed \emph{away} from the vertex. For this paper it is convenient to fix a formal definition of a walk (and related notions like circuits and cycles) using single transitions. A \emph{closed walk} $W$ is a sequence $((h_1,h_2),\dots,(h_{n-1},h_n))$ of directed single transitions, where $\{h_n,h_1\}$ is an edge and $\{h_i,h_{i+1}\}$ is an edge if $i \in \{1,\dots,n-1\}$ is even. We consider closed walks modulo cyclic shifts, i.e., we assume that closed walks have no distinguished starting vertex. We say that $W$ traverses $e$ if the half-edges of $e$ appear in some directed single transitions of $W$. A nonempty closed walk $W$ is an \emph{oriented circuit} if each half-edge appears at most once in $W$. Thus an oriented circuit visits each edge at most once, but there may be vertex repetitions. We (may) conveniently represent an oriented circuit $((h_1,h_2),\dots,(h_{n-1},h_n))$ as the nonempty set $\{(h_1,h_2),\dots,(h_{n-1},h_n)\}$ since we assume no distinguished starting vertex. An \emph{oriented cycle} is a set $C$ of oriented circuits such that every edge is traversed by at most one oriented circuit from $C$.

The notions of a \emph{cycle} and \emph{circuit} capture the notions of an oriented cycle and an oriented circuit where we additionally forget the orientation, i.e., each directed single transition is replaced by its corresponding single transition. Note that, just like an oriented circuit, a circuit may visit a vertex more than once.

\subsection{Cycle bases}

Let us recall the well-known notion of an incidence matrix.
\begin{definition}
The \emph{incidence matrix} of a directed graph $D$ is the $V(D) \times E(D)$-matrix $M$ over $\mathbb{Q}$ where entry $M_{v,e}$ with $v \in V(D)$ and $e \in E(D)$ is $1$ if $v$ is incident to the tail but not the head of $e$, $-1$ if $v$ is incident to the head but not the tail of $e$, and $0$ otherwise.
\end{definition}
The \emph{cycle space} of $D$ is the right nullspace of its incidence matrix.

Let $G$ be a graph and $D$ be a directed version of $G$. For a closed walk $W$ of $G$, we let $\eCnt{D}{W} \in \mathbb{Z}^{E(D)}$ be obtained from the zero vector by tallying $+1$ ($-1$, resp.)\ for the entry with index $e$ each time $W$ traverses $e$ along (against, resp.)\ the direction of $D$. We call $\eCnt{D}{W}$ the \emph{incidence vector} of $W$ in $D$. For a set $S$ of closed walks, we write $\eCnt{D}{S} = \{\eCnt{D}{W} \mid W \in S\}$.

The cycle space of $D$ is equal to $\spn_\mathbb{Q}(\eCnt{D}{S})$, where $S$ is the set of oriented circuits of $G$. Note that we can equivalently define the notion of cycle space in terms of closed walks or oriented cycles instead of oriented circuits. The following elementary property of cycle spaces will be useful later.

\begin{lemma}\label{lem:sum_oriented_cycles}
Let $D$ be a directed version of a graph $G$ and $s \in \mathbb{Q}^{E(D)}$. There is an oriented cycle $C$ of $G$ with $s=\eCnt{D}{C}$ if and only if $s$ is an element of the cycle space of $D$ and every entry of $s$ is in $\{-1,0,1\}$.
\end{lemma}
\begin{proof}
If $s=\eCnt{D}{C}$ for some oriented cycle $C$, then certainly $s$ is an element of the cycle space of $D$ and every entry of $s$ is in $\{-1,0,1\}$.

The converse is verified by induction on the number of nonzero entries of $s$. For simplicity, we reverse the direction of every edge of $D$ whose corresponding entry in $s$ is $-1$, and proceed with the assumption that every entry of $s$ is in $\{0,1\}$. If $s$ is the zero vector, then $s=\eCnt{D}{\emptyset}$. Otherwise, let $e_1$ be an edge with $s(e_1) = 1$. Let $h_1$ be the head of $e_1$. Since $s$ is in the right nullspace of the incidence matrix of $D$, the vertex $v_1$ incident to $h_1$ is also incident to a half-edge $h_2$ that is the tail of an edge $e_2 = (h_2,h_3) \in E(D)$ and has $s(e_2)=1$. We observe that $(h_1,h_2)$ is a directed single transition. Continuing in this fashion with the head $h_3$ of $e_2$, we obtain an oriented circuit $C_1 = \{ (h_1,h_2), (h_3,h_4), \ldots, (h_{n-1},h_n) \}$. Now, $s-\eCnt{D}{C_1}$ is an element of the cycle space of $D$ having the same entries as $s$ for the edges not traversed by $C_1$, and zero entries for the edges traversed by $C_1$. By the inductive hypothesis, the assertion of this lemma applies to $s - \eCnt{D}{C_1}$, so it also applies to $s = (s-\eCnt{D}{C_1})+\eCnt{D}{C_1}$.
\end{proof}

A \emph{cycle basis} of $D$ is a set $B$ of closed walks of $G$ such that $\eCnt{D}{B}$ is of cardinality $|B|$ and forms a basis of the cycle space of $D$. 

\begin{remark}\label{rem:bases_vertex_rep}
We remark that the notion of cycle basis as defined here in terms of closed walks is more general than usual in the literature. In particular, we allow vertex and edge repetitions (i.e., a closed walk may visit a vertex or edge more than once). We need to allow vertex repetitions when we consider oriented circuits induced by Eulerian circuits in Section~\ref{sec:cycle_basis_Euler}.
\end{remark}

Note that if $D_1$ and $D_2$ are directed versions of $G$, then $B$ is a cycle basis of $D_1$ if and only if $B$ is a cycle basis of $D_2$. Therefore, we (may) speak of a \emph{cycle basis} of $G$. Similarly, a \emph{cycle spanning set} of $G$ is a set $I$ of closed walks of $G$ such that there is a subset $B$ of $I$ that is a cycle basis of $G$. Since a maximal forest of $G$ has $|V(G)| - c(G)$ edges, for any cycle basis $B$ of $G$, we have $|B| = |E(G)| - (|V(G)| - c(G))$.

We say that a cycle spanning set $B$ of $D$ is \emph{integral} if every closed walk $W$ of $D$ has $\eCnt{D}{W} \in \spn_\mathbb{Z}(\eCnt{D}{B})$. That is, a cycle spanning set $B$ is integral if for each closed walk $W$, we have that $\eCnt{D}{W}$ is a linear combination of elements from $\eCnt{D}{B}$ using integer coefficients, i.e., $\eCnt{D}{W} = \sum_{W' \in B} \lambda_{W'} \eCnt{D}{C'}$ where $\lambda_{W'} \in \mathbb{Z}$ for all $W' \in B$. An \emph{integral cycle basis} is an integral cycle spanning set that is also a cycle basis.

Note that the notion of integral cycle basis is also independent of the chosen directed version $D$ of $G$. Hence we (may) speak of an integral cycle basis (or spanning set) of $G$.

\section{Cycle bases of touch-graphs}\label{sec:touch_graphs}

Recall that a graph is called \emph{$k$-regular} if every vertex is incident to exactly $k$ half-edges. Since a circle graph is uniquely determined by fixing Euler circuits for each connected component of a 4-regular graph (see the definition of a circle graph in Section~\ref{sec:cycle_matrices} below), the theory of circle graphs is intimately connected to the theory of 4-regular graphs. 

In this section we recall the definition of the touch-graph of a circuit partition of a 4-regular graph. We also introduce the notion of an edge-transition incidence matrix, which is closely related to the notion of an incidence matrix and is central in the proof of the main result of this paper (Theorem~\ref{thm:main_result}).

A \emph{circuit partition} $P$ of a 4-regular graph $F$ is a cycle such that every half-edge of $F$ occurs in exactly one circuit of $P$. So, informally, $P$ partitions the edges of $F$ into circuits. Circuit partitions can be described in terms of transitions. If $F$ is 4-regular then for a vertex $v \in V(F)$, a \emph{transition} at $v$ is a partition of the set of half-edges incident to $v$ in pairs; equivalently, it is a pair of disjoint single transitions at $v$. The set of transitions of $F$ is denoted $\mathfrak{T}(F)$. A \emph{transversal} of $\mathfrak{T}(F)$ contains exactly one transition for each vertex of $F$.

For a circuit partition $P$, denote by $\tau(P)$ the transversal of $\mathfrak{T}(F)$ that includes the transitions $t$ corresponding to $P$ (i.e., each single transition of $t$ is in a circuit of $P$). Conversely, each transversal $T$ of $\mathfrak{T}(F)$ uniquely determines a circuit partition $P$ with $\tau(P)=T$.

Each transition $t \in \tau(P)$ corresponds to an edge in a graph called the touch-graph of $P$  \cite{bouchet1988}. The touch-graph encodes the incidences between the circuits of $P$ and the single transitions of $F$.

\begin{definition}
Let $P$ be a circuit partition of a $4$-regular graph $F$. Then the \emph{touch-graph} of $P$, denoted by $\Tch(P)$, is the graph $(P,\xi(P),\tau(P),\epsilon)$, where $\xi(P) := \bigcup P$ is the set of single transitions corresponding to $P$ and $\epsilon$ maps every $s \in \xi(P)$ to the $C \in P$ such that $s \in C$.
\end{definition}

According to the definition, the edges of $\Tch(P)$ correspond to elements of $\tau(P)$. As $\tau(P)$ has one element for each vertex of $F$, the edges of $\Tch(P)$ also correspond to vertices of $F$. Therefore the touch-graphs of the circuit partitions of $F$ are all related to each other through bijections of their edges.

\newcommand{\Ccross}[1]{\path[draw,color=blue]
    (#1.45) to (#1.225)
    (#1.135) to (#1.315);
}
\newcommand{\Clr}[1]{\path[draw,color=blue]
    (#1.45) [out=225,in=315] to (#1.135)
    (#1.225) [out=45,in=135] to (#1.315);
}
\newcommand{\Ctd}[1]{\path[draw,color=blue]
    (#1.135) [out=315,in=45] to (#1.225)
    (#1.45) [out=225,in=135] to (#1.315);
}
\newcommand{\fourreg}{
\tikzstyle{vert}=[circle,minimum size=14,inner sep=0pt,draw]
\node[vert,label={above:$c$}] (pc) at (0,2.7) {};
\node[vert,label={left:$d$}] (pd) at (0,1.2) {};
\node[vert,label={left:$a$}] (pa) at (-1.2,0) {};
\node[vert,label={right:$b$}] (pb) at (1.2,0) {};
\draw[edg]
(pa) edge [out=315,in=225,above] node {$e_1$} (pb)
(pb) edge [out=45,in=45,right] node {$e_2$} (pc)
(pc) edge [out=225,in=135,left] node {$e_3$} (pd)
(pd) edge [out=315,in=135,right] node {$e_4$} (pb)
(pb) edge [out=315,in=225,above,looseness=2] node {$e_5$} (pa)
(pa) edge [out=135,in=135,left] node {$e_6$} (pc)
(pc) edge [out=315,in=45,right] node {$e_7$} (pd)
(pd) edge [out=225,in=45,right,near end] node {$e_8$} (pa)
;
}
\begin{figure}
\begin{center}
\begin{tikzpicture}[scale=0.9]
\tikzset{edg/.style={}}
\fourreg
\end{tikzpicture}
\begin{tikzpicture}[scale=0.9]
\tikzset{edg/.style={}}
\fourreg
\Ccross{pc}
\Clr{pd}
\Ccross{pa}
\Ccross{pb}
\end{tikzpicture}
\end{center}
\caption{A $4$-regular graph $F$ (left) and a circuit partition $P$ of $F$ (right).}
\label{fig:4reg_graph}
\end{figure}

\begin{figure}
\begin{center}
\begin{tikzpicture}
\tikzstyle{vert}=[circle,minimum size=14,inner sep=0pt,draw]
\node[vert] (x) at (0,0) {};
\node[vert] (y) at (2,0) {};
\path[draw]
(x) to [bend left=50,above] node {$a$} (y)
(x) to [above] node {$b$} (y)
(x) to [bend right,below] node {$d$} (y)
(x) to [loop left,in=135,out=225,looseness=8] node {$c$} (x)
;
\end{tikzpicture}
\end{center}
\caption{The graph $\Tch(P)$ with $P$ from Figure~\ref{fig:4reg_graph}.}
\label{fig:tch}
\end{figure}

\begin{example}\label{ex:4reg_graph}
Consider the $4$-regular graph $F$ on the left-hand side of Figure~\ref{fig:4reg_graph}. We use $F$ as a running example. The right-hand side of this figure represents a circuit partition $P$ of $F$ by depicting $\tau(P)$ in blue. We notice that $|P|=2$, so $\Tch(P)$ has two vertices. Moreover, the single transitions of the transition $t_c$ at $c$ in $\tau(P)$ belong to a common circuit of $P$, so $t_c$ is a loop in $\Tch(P)$. The single transitions of the transitions in $\tau(P)$ at $a$, $b$, and $d$ belong to different circuits of $P$, so these transitions are non-loop edges in $\Tch(P)$. The graph $\Tch(P)$ is depicted in Figure~\ref{fig:tch}, where, for notational convenience, instead of the edge identities (i.e., transitions) the figure gives the vertices at which these transitions reside.
\end{example}

It is useful to have a notation for transitions with respect to Eulerian circuits. For an Eulerian circuit $C$ of a connected $4$-regular graph $F$ and a vertex $v$ of $F$, we denote by $\phi_C(v)$ the transition at $v$ that is included in $\tau(C)$. Suppose that the directed single transitions $(h_1,h_2)$ and $(h_1',h_2')$ both appear at $v$ in one of the orientations of $C$. Then we denote by $\chi_C(v)$ the transition $\{\{h_1,h_2'\},\{h_1',h_2\}\}$, and we denote by $\psi_C(v)$ the transition $\{\{h_1,h_1'\},\{h_2,h_2'\}\}$. Note that the notions $\phi_C(v)$, $\chi_C(v)$, and $\psi_C(v)$ are independent of the chosen orientation of $C$. For a  4-regular graph $F$, an \emph{Euler system} is a set containing, for each connected component $F'$ of $F$, exactly one Eulerian circuit of $F'$. Given an Euler system $C$ of $F$, we define, for vertices $v$ of $F$, $\phi_C(v) := \phi_{C'}(v)$, $\chi_C(v) := \chi_{C'}(v)$, and $\psi_C(v) := \psi_{C'}(v)$, where $C' \in C$ is the Euler circuit of the connected component containing $v$.

\begin{figure}
\begin{center}
\begin{tikzpicture}[scale=0.9]
\tikzset{edg/.style={}}
\fourreg
\Ctd{pa}
\Ccross{pb}
\Ccross{pc}
\Ccross{pd}
\end{tikzpicture}
\begin{tikzpicture}[scale=0.9]
\tikzset{edg/.style={>=latex,->}}
\fourreg
\Ctd{pa}
\Ccross{pb}
\Ccross{pc}
\Ccross{pd}
\end{tikzpicture}
\end{center}
\caption{An Eulerian circuit $C$ of $F$ (left) and an orientation $D$ of $C$ (right).}
\label{fig:Euler}
\end{figure}

\begin{example}\label{ex:P_in_terms_of_C}
Consider again the $4$-regular graph $F$ of Example~\ref{ex:4reg_graph}. An Eulerian circuit $C$ of $F$ is depicted on the left-hand side of Figure~\ref{fig:Euler} and an orientation of $C$ is depicted on the right-hand side of this figure. The transversal corresponding to circuit partition $P$ in Figure~\ref{fig:4reg_graph} is $\{\psi_C(a),\phi_C(b),\phi_C(c),\psi_C(d)\}$.
\end{example}

Let $F$ be a $4$-regular graph. A \emph{transitional orientation} $\mathbf{o}$ of $F$ is a function that assigns to each transition $t \in \mathfrak{T}(F)$ one of its two single transitions $\mathbf{o}(t) \in t$. We now introduce a notion that is somewhat similar to the notion of an incidence matrix and is central to our proof of Theorem~\ref{thm:main_result} (cf.\ Corollary~\ref{cor:trans_totunimod} below).

\begin{definition}
Let $F$ be a $4$-regular graph and let $D$ be a directed version of $F$. Let $\mathbf{o}$ be a transitional orientation of $F$.

The \emph{edge-transition incidence matrix} of $D$ with respect to $\mathbf{o}$, denoted by $\edgetrans{D}{\mathbf{o}}$, is the $E(F) \times \mathfrak{T}(F)$-matrix over $\mathbb{Q}$ where, for each $e \in E(F)$ and each $t \in \mathfrak{T}(F)$, its entry indexed by $(e,t)$ is 
\[
\begin{cases}
1 & e \cap \mathbf{o}(t) = \{h\} \text{ and $h$ is the tail of $e$ in $D$,}\cr
-1 & e \cap \mathbf{o}(t) = \{h\} \text{ and $h$ is the head of $e$ in $D$,}\cr
0 & \text{otherwise.}
\end{cases}
\]
\end{definition}

Notice that the column of $\edgetrans{D}{\mathbf{o}}$ corresponding to a transition $t$ has two nonzero entries unless $\mathbf{o}(t)$ is a loop, in which case the $t$-column is $0$. Also, if $D_1$ and $D_2$ are directed versions of $F$, then $\edgetrans{D_2}{\mathbf{o}}$ can be obtained from $\edgetrans{D_1}{\mathbf{o}}$ by multiplying the $e$-rows where $e$ has different orientations in $D_1$ and $D_2$ by $-1$. In contrast, if $\mathbf{o}_1$ and $\mathbf{o}_2$ are transitional orientations, then it is not so easy to describe the connection between $\edgetrans{D}{\mathbf{o}_1}$ and $\edgetrans{D}{\mathbf{o}_2}$. For instance, they may have different numbers of zero columns, and different numbers of zero rows.

Moreover, notice the strong similarity between the incidence matrix of $D$ and $\edgetrans{D}{\mathbf{o}}$. Very roughly (in particular, assuming no loops), if $t$ is a transition at $v$ then the column with index $t$ of $\edgetrans{D}{\mathbf{o}}$ is obtained from the row with index $v$ of the incidence matrix of $D$ by setting two of the four nonzero entries to zero. Here $\mathbf{o}$ determines which two entries are set to zero.

Finally, we notice that $\edgetrans{D}{\mathbf{o}}$ is the product of an $E(F) \times H(F)$-matrix $H_1$ and an $H(F) \times \mathfrak{T}(F)$-matrix $H_2$, where (1) for $e \in E(F)$ and $h \in H(F)$, the entry of $H_1$ indexed by $(e,h)$ is $1$ if $h$ is the tail of $e$ in $D$, $-1$ if $h$ is the head of $e$ in $D$, and $0$ otherwise, and (2) for $h \in H(F)$ and $t \in \mathfrak{T}(F)$ the entry of $H_2$ indexed by $(h,t)$ is $1$ if $h \in \mathbf{o}(t)$ and $0$ otherwise.

\begin{example}\label{ex:eti}
Consider again the $4$-regular graph $F$ and the Euler system $C$ from the left-hand side of Figure~\ref{fig:Euler}. Consider the directed version $D$ of $F$ induced by the orientation of $C$ as depicted on the right-hand side of Figure~\ref{fig:Euler}. Then the transpose of the incidence matrix of $D$ is as follows:
\[
\bordermatrix{
~   &  a &  b &  c &  d \cr
e_1 &  1 & -1 &  0 &  0 \cr
e_2 &  0 &  1 & -1 &  0 \cr
e_3 &  0 &  0 &  1 & -1 \cr
e_4 &  0 & -1 &  0 &  1 \cr
e_5 & -1 &  1 &  0 &  0 \cr
e_6 &  1 &  0 & -1 &  0 \cr
e_7 &  0 &  0 &  1 & -1 \cr
e_8 & -1 &  0 &  0 &  1
}.
\]
Let $\mathbf{o}$ be the transitional orientation that assigns to a transition at $v \in V(F)$ the single transition that does not contain any of the heads of edges $e_4$, $e_5$, $e_6$, and $e_7$ (for vertices $b$, $a$, $c$, and $d$, respectively). This particular $\mathbf{o}$ is chosen in view of Example~\ref{ex:IAS_from_Euler} towards the end of the paper. Now $\edgetrans{D}{\mathbf{o}}$ is as follows:
\[
\scalebox{0.9}{
{\let\quad\thinspace
\bordermatrix{
~ & 
\phi_C(a) & \phi_C(b) & \phi_C(c) & \phi_C(d) &  
\chi_C(a) & \chi_C(b) & \chi_C(c) & \chi_C(d) & 
\psi_C(a) & \psi_C(b) & \psi_C(c) & \psi_C(d) \cr
e_1 &  1 & -1 &  0 &  0 &  0 & -1 &  0 &  0 &  1 &  0 &  0 &  0 \cr
e_2 &  0 &  1 & -1 &  0 &  0 &  0 & -1 &  0 &  0 &  1 &  0 &  0 \cr
e_3 &  0 &  0 &  1 & -1 &  0 &  0 &  0 & -1 &  0 &  0 &  1 &  0 \cr
e_4 &  0 &  0 &  0 &  1 &  0 &  0 &  0 &  0 &  0 &  0 &  0 &  1 \cr
e_5 &  0 &  0 &  0 &  0 &  0 &  1 &  0 &  0 &  0 &  1 &  0 &  0 \cr
e_6 &  0 &  0 &  0 &  0 &  1 &  0 &  0 &  0 &  1 &  0 &  0 &  0 \cr
e_7 &  0 &  0 &  0 &  0 &  0 &  0 &  1 &  0 &  0 &  0 &  1 &  0 \cr
e_8 & -1 &  0 &  0 &  0 & -1 &  0 &  0 &  1 &  0 &  0 &  0 &  1 
}}}.
\]
\end{example}

For a finite set $X$, an $X$-vector $v$, and $Y \subseteq X$, we let $v|_Y$ denote the $Y$-vector obtained from $v$ by restricting to the entries of $Y$. Similarly, for a $W \times X$-matrix $A$ and $Y \subseteq X$, we let $A|_Y$ denote the $W \times Y$-matrix obtained from $A$ by restricting to the columns of $Y$. 

A transitional orientation $\mathbf{o}$ is used to simultaneously fix directions of the edges of $\Tch(P)$ for all circuit partitions $P$ of a $4$-regular graph $F$. For a circuit partition $P$ of $F$, we denote by $\dTch{P}{\mathbf{o}}$ the directed version of $\Tch(P)$ where each edge $t \in \tau(P)$ is directed from the $p \in P$ containing the single transition of $t$ distinct from $\mathbf{o}(t)$ towards the $p' \in P$ containing $\mathbf{o}(t)$.

\begin{figure}
\begin{center}
\begin{tikzpicture}
\tikzset{edg/.style={>=latex,->}}
\tikzstyle{vert}=[circle,minimum size=14,inner sep=0pt,draw]
\node[vert] (x) at (0,0) {};
\node[vert] (y) at (2,0) {};
\path[edg]
(y) edge [bend right=50,above] node {$a$} (x)
(y) edge [above] node {$b$} (x)
(x) edge [bend right,below] node {$d$} (y)
(x) edge [loop left,in=135,out=225,looseness=8] node {$c$} (x)
;
\end{tikzpicture}
\end{center}
\caption{The graph $\dTch{P}{\mathbf{o}}$ from Example~\ref{ex:dir_tch}.}
\label{fig:dir_tch}
\end{figure}

\begin{example}\label{ex:dir_tch}
Consider again $F$, $P$, $D$, and $\mathbf{o}$ from the running example; see Figure~\ref{fig:4reg_graph} for $F$ and $P$, Figure~\ref{fig:Euler} for $D$, and Example~\ref{ex:eti} for $\mathbf{o}$. Then $\dTch{P}{\mathbf{o}}$ is given in Figure~\ref{fig:dir_tch}, where the loop $(s_1, s_2)$ corresponding to $c$ is such that $s_1$ contains the head of $e_6$ and the tail of $e_7$ and $s_2$ contains the head of $e_2$ and the tail of $e_3$.
\end{example}

For a $4$-regular graph $F$ and a circuit partition $P$, a closed walk $W$ of $F$ determines a closed walk of $\Tch(P)$, denoted by $\pi_P(W)$, as follows (see also \cite{Traldi/SignedInterlace}): a visit of $W$ at a vertex $v$ that traverses a single transition not in the transition $t$ of $P$ at $v$ corresponds to walking the edge $t$ of $\Tch(P)$ (in this way ``jumping'' from one position in a circuit of $P$ to another position in a circuit of $P$, possibly within the same circuit), while traversing a single transition of $t$ corresponds to ``staying put'' in the vertex of $\Tch(P)$ (i.e., circuit of $P$) that contains this single transition. More precisely, consider the sequence $W'$ obtained from $W$ by replacing every directed single transition $(h,h')$ in $W$ by the tuple $(s,s')$, where $s$ and $s'$ are single transitions of circuits of $P$ such that $h \in s$ and $h' \in s'$ (note that $s$ and $s'$ are unique with this property). Now, $\pi_P(W)$ is obtained from $W'$ by removing all tuples of the form $(s,s)$. Note that, for every tuple $(s,s')$ of $\pi_P(W)$, we have $\{s,s'\} \in \tau(P) = E(\Tch(P))$ (since $s \neq s'$) and that $\pi_P(W)$ is indeed a closed walk of $\Tch(P)$. 

\begin{example}
Consider again $F$ and $P$ from Figure~\ref{fig:4reg_graph}. Consider the closed walk $W$ of $F$ that traverses exactly once the edge $e_1$ in the direction of $b$, followed by the edge $e_4$, and finally the edge $e_8$. Then $\pi_P(W)$ is the closed walk that traverses the edge corresponding to $b$ in the direction of the vertex without the loop, followed by traversing the edge corresponding to $a$. Note that the visit of $W$ at vertex $d$ does not correspond to any edge traversal in $\pi_P(W)$.
\end{example}

Similarly as for $\eCnt{D}{S}$, we write $\pi_P(S) = \{\pi_P(W) \mid W \in S\}$ for a set $S$ of closed walks. In fact, we regard here $\pi_P(S)$ as a multiset to ensure $|\pi_P(S)| = |S|$, which will be important when we turn to matrices in Section~\ref{sec:cycle_matrices}. The definition of $\pi_P$ carries over in the natural way to (sets of) ``unoriented'' closed walks (i.e., closed walks where we forget the orientation). By the definition of $\pi_P$, one observes that, for every closed walk $W$ of $\Tch(P)$, there is a closed walk $W'$ such that $\pi_P(W') = W$.

\begin{theorem}\label{thm:pi_on_vectors_full}
Let $F$ be a $4$-regular graph and let $D$ be a directed version of $F$. Let $\mathbf{o}$ be a transitional orientation of $F$. Let $P$ be a circuit partition of $F$. 

We have 
\[
\eCnt{D}{W} \cdot \edgetrans{D}{\mathbf{o}}|_{\tau(P)} = \eCnt{\dTch{P}{\mathbf{o}}}{\pi_P(W)}
\]
for all closed walks $W$ of $F$. In this equality, the $\eCnt{\cdot}{\cdot}$ vectors are interpreted as row vectors.
\end{theorem}
\begin{proof}
Let $W$ be a closed walk of $F$ and let $t_v \in \tau(P)$ be a transition at a vertex $v \in V(F)$. If $W$ does not visit $v$, then the $t_v$-entries of $\eCnt{D}{W} \cdot \edgetrans{D}{\mathbf{o}}|_{\tau(P)}$ and $\eCnt{\dTch{P}{\mathbf{o}}}{\pi_P(W)}$ are both zero. Let us consider a visit of $v$ by $W$. Assume that this visit traverses the single transition $s = \{h_1,h_2\}$ in the direction $(h_1,h_2)$, i.e., this visit arrives at $v$ via half-edge $h_1$ and leaves $v$ via half-edge $h_2$. Let $e_1$ and $e_2$ be the edges of $F$ corresponding to $h_1$ and $h_2$, respectively. The decomposition of $\edgetrans{D}{\mathbf{o}}$ as the product of $H_1$ and $H_2$ (see above Example~\ref{ex:eti}) implies (by restricting the columns of $H_1$ and the rows of $H_2$ to $s$) that the contribution of this visit to the $t_v$-entry of $\eCnt{D}{W} \cdot \edgetrans{D}{\mathbf{o}}|_{\tau(P)}$ is $-1$ for the ``incoming'' edge $e_1$ if $h_1 \in \mathbf{o}(t_v)$ and $0$ otherwise, and $1$ for the ``outgoing'' edge $e_2$ if $h_2 \in \mathbf{o}(t_v)$ and $0$ otherwise. (If $e_1 = e_2$ and $\mathbf{o}(t_v) = \{h_1,h_2\}$, then the total contribution of this visit is $0$.)

If $s \in t_v$, then the contribution of this visit is $0$ --- as required. Assume now that $s \notin t_v$. Then exactly one of $h_1$ and $h_2$ is in $\mathbf{o}(t_v)$.

If $h_1 \in \mathbf{o}(t_v)$, then the contribution of this visit to the $t_v$-entry of $\eCnt{D}{W} \cdot \edgetrans{D}{\mathbf{o}}|_{\tau(P)}$ is $-1$, which corresponds to traversing the edge $t_v$ of $\Tch(P)$ against the direction of $\dTch{P}{\mathbf{o}}$. Since the direction of $t_v$ in $\dTch{P}{\mathbf{o}}$ is from the unique single transition $s' \in t_v \setminus \{\mathbf{o}(t_v)\}$ to $\mathbf{o}(t_v)$, the corresponding edge visit in $\pi_P(W)$ is indeed against the direction of $\dTch{P}{\mathbf{o}}$.

Similarly, if $h_2 \in \mathbf{o}(t_v)$, then the contribution of this visit to the $t_v$-entry of $\eCnt{D}{W} \cdot \edgetrans{D}{\mathbf{o}}|_{\tau(P)}$ is $1$, which corresponds to traversing the edge $t_v$ of $\Tch(P)$ along the direction of $\dTch{P}{\mathbf{o}}$. Since the direction of $t_v$ in $\dTch{P}{\mathbf{o}}$ is from $s'$ to $\mathbf{o}(t_v)$, the corresponding edge visit in $\pi_P(W)$ is indeed along the direction of $\dTch{P}{\mathbf{o}}$.
\end{proof}

\begin{figure}
\[
\begin{tikzcd}[column sep=huge]
W \arrow{r}{\eCnt{D}{\cdot}} \arrow[swap]{d}{\pi_P} & \eCnt{D}{W} \arrow{d}{\edgetrans{D}{\mathbf{o}}|_{\tau(P)}} \\
\pi_P(W) \arrow{r}{\eCnt{\dTch{P}{\mathbf{o}}}{\cdot}} & \eCnt{\dTch{P}{\mathbf{o}}}{\pi_P(W)}
\end{tikzcd}
\]
\caption{Commutative diagram for closed walks $W$ of $4$-regular graphs.}\label{fig:comm_diagram}
\end{figure}

Therefore $\edgetrans{D}{\mathbf{o}}$ corresponds to a linear transformation sending incidence vectors of $D$ to incidence vectors of $\Tch_\mathbf{o}(P)$ in a way compatible with $\pi_P$; see Figure~\ref{fig:comm_diagram}. Also note that the left-hand side of the equality of Theorem~\ref{thm:pi_on_vectors_full} depends on $D$, but the right-hand side of this equality does not.

\begin{figure}
    \begin{center}
    \begin{tikzpicture}
        \tikzstyle{vert}=[circle,minimum size=14,inner sep=0pt,draw]
        \tikzset{edg/.style={>=latex,->}}

\node[vert] (x) at (0,0) {};
        \draw[edg]
        (x) edge [loop left,in=135,out=225,looseness=12] node {$e_2$} node[above,pos=0.95] {$h_3$} node[below,pos=0.05] {$h_2$} (x)
        (x) edge [loop right,in=45,out=-45,looseness=12] node {$e_1$} node[above,pos=0.95] {$h_1$} node[below,pos=0.05] {$h_4$} (x)
        ;
        \Ccross{x}
        \node [below=1cm, align=flush center] at (x)
        {
            (\emph{i})
        };
    \end{tikzpicture}\\
    \vspace{0.3cm}
    \begin{tikzpicture}[scale=0.9, every node/.style={scale=0.9}]
        \tikzstyle{vert}=[circle,minimum size=14,inner sep=0pt,draw]
        \tikzset{edg/.style={>=latex,->}}
\node[vert] (x) at (0,0) {};
        \draw[edg]
        (x) edge [loop above,in=45,out=135,looseness=12] node {$t_1$} node[right,pos=0.95] {$\{h_1,h_2\}$} node[left,pos=0.05] {$\{h_3,h_4\}$} (x)
        ;
        \node [below=0.7cm, align=flush center] at (x)
        {
            (\textit{ii})
        };

\node[vert] (x) at (4,0) {};
        \draw[edg]
        (x) edge [loop above,in=45,out=135,looseness=12] node {$t_2$} node[right,pos=0.95] {$\{h_1,h_3\}$} node[left,pos=0.05] {$\{h_2,h_4\}$} (x)
        ;
        \node [below=0.7cm, align=flush center] at (x)
        {
            (\textit{iii})
        };

\node[vert] (x) at (7,0.5) {};
        \node[vert] (y) at (11,0.5) {};
        \path[edg]
        (x) edge [above] node {$t_3$} node[below,pos=0.8] {$\{h_1,h_4\}$} node[below,pos=0.2] {$\{h_2,h_3\}$} (y)
        ;
        \node [below=0.7cm, align=flush center] at (9,0)
        {
            (\textit{iv})
        };
    \end{tikzpicture}
    \end{center}
    \caption{Some graphs, with edge and half-edge identities, considered in Example~\ref{ex:one_vert}: (\textit{i}) the directed version $D$ of $F$, where the directed single transitions of $W$ are depicted inside the vertex, and the graphs (\textit{ii}) $\Tch_\mathbf{o}(P_1)$, (\textit{iii}) $\Tch_\mathbf{o}(P_2)$, and (\textit{iv}) $\Tch_\mathbf{o}(P_3)$.}
    \label{fig:one_vert}
\end{figure}

To illustrate Theorem~\ref{thm:pi_on_vectors_full}, we give two examples: one where $F$ is the 4-regular graph with only one vertex to cover the case of loops, and one where we continue the running example.
\begin{example}\label{ex:one_vert}
Suppose $F$ is the 4-regular graph with only one vertex and two edges $e_1 = \{h_1, h_4\}$ and $e_2 = \{h_2, h_3\}$ (both loops). Then $W = ((h_1, h_2), (h_3, h_4))$ is a closed walk of $F$. Suppose that $D$ is the directed version of $F$ consistent with the orientation of $W$, depicted in Figure~\ref{fig:one_vert}(\textit{i}). Let $\mathbf{o}$ be the transitional orientation of $F$ with $h_1 \in \mathbf{o}(t)$ for every transition $t$. If the transitions of $F$ are $t_1, t_2, t_3$ with $\mathbf{o}(t_i) = \{h_1, h_{i+1} \}$, then
\[
\edgetrans{D}{\mathbf{o}} =
\bordermatrix{
~ & t_1 & t_2 & t_3 \cr
e_1 & -1 & -1 & 0 \cr
e_2 & 1 & -1 & 0
}.
\]
Now $\sigma(D,W) = \begin{pmatrix} 1 & 1 \end{pmatrix}$. For $i \in \{1, 2, 3\}$, let $P_i$ be the circuit partition of $F$ determined by the transition $t_i$. Then Theorem~\ref{thm:pi_on_vectors_full} is satisfied because $\pi_{P_1}(W)$ is the empty closed walk of $\Tch_\mathbf{o}(P_1)$; $\pi_{P_2}(W)$ is the closed walk that traverses the loop edge of $\Tch_\mathbf{o}(P_2)$ twice, against its direction; and $\pi_{P_3}(W)$ is a closed walk that traverses the non-loop edge of $\Tch_\mathbf{o}(P_3)$ twice, in opposite directions. See (\textit{ii})-(\textit{iv}) of Figure~\ref{fig:one_vert} for depictions of these directed touch graphs.
\end{example}

\begin{example}
Consider again $F$, $P$, $C$, $D$, and $\mathbf{o}$ from the running example; see Figure~\ref{fig:4reg_graph} for $F$ and $P$, Figure~\ref{fig:Euler} for $C$ and $D$, and Example~\ref{ex:eti} for $\mathbf{o}$. Consider again the closed walk $W$ of $F$ that traverses exactly once the edge $e_1$ along the direction of $D$, followed by the edge $e_4$, and finally the edge $e_8$. We have
\[
\eCnt{D}{W} =
\bordermatrix{
& e_1 & e_2 & e_3 & e_4 & e_5 & e_6 & e_7 & e_8 \cr
&   1 &   0 &   0 &  -1 &   0 &   0 &   0 &   1
},
\]
and therefore $\eCnt{D}{W} \cdot \edgetrans{D}{\mathbf{o}}$ is equal to
\[
\scalebox{0.9}{
{\let\quad\thinspace
\bordermatrix{
& 
\phi_C(a) & \phi_C(b) & \phi_C(c) & \phi_C(d) &  
\chi_C(a) & \chi_C(b) & \chi_C(c) & \chi_C(d) & 
\psi_C(a) & \psi_C(b) & \psi_C(c) & \psi_C(d) \cr
&  0 & -1 &  0 & -1 & -1 & -1 &  0 &  1 &  1 &  0 &  0 &  0
}}}.
\]
By Example~\ref{ex:P_in_terms_of_C} and Theorem~\ref{thm:pi_on_vectors_full},
\[
\eCnt{\dTch{P}{\mathbf{o}}}{\pi_P(W)} = 
\bordermatrix{
& 
\psi_C(a) & \phi_C(b) & \phi_C(c) & \psi_C(d) \cr
&  1 & -1 &  0 &  0
}.
\]
Since $F$ in this example has no loops, we (can) unambiguously denote the half edge of an edge $e$ incident to a vertex $v$ by $h_{e,v}$. Then $W = ((h_{e_1,b}, h_{e_4,b}),\allowbreak (h_{e_4,d}, h_{e_8,d}),\allowbreak (h_{e_8,a}, h_{e_1,a}))$.

Note that $\dTch{P}{\mathbf{o}}$ is obtained from $\Tch(P)$ (depicted in Figure~\ref{fig:tch}) by directing $\{h_{e_5,a}, h_{e_8,a}\}$ to $\{h_{e_1,a}, h_{e_6,a}\}$, $\{h_{e_4,b}, h_{e_5,b}\}$ to $\{h_{e_1,b}, h_{e_2,b}\}$, $\{h_{e_6,c}, h_{e_7,c}\}$ to $\{h_{e_2,c}, h_{e_3,c}\}$, and $\{h_{e_3,d}, h_{e_7,d}\}$ to $\{h_{e_4,d}, h_{e_8,d}\}$.

We have
\[
\pi_P(W) = ((\{h_{e_1,b}, h_{e_2,b}\}, \{h_{e_4,b}, h_{e_5,b}\}), (\{h_{e_5,a}, h_{e_8,a}\}, \{h_{e_1,a}, h_{e_6,a}\})),
\]
and so $\pi_P(W)$ indeed traverses the edge corresponding to $b$ against the direction of $\dTch{P}{\mathbf{o}}$ and the edge corresponding to $a$ along the direction of $\dTch{P}{\mathbf{o}}$.
\end{example}

\begin{corollary}\label{cor:transfer_weights}
Let $F$ be a $4$-regular graph, let $D$ be a directed version of $F$, let $W$ be a closed walk of $F$, let $P$ be a circuit partition of $F$, and let $D'$ be a directed version of $\Tch(P)$. If 
\[
\eCnt{D}{W} = \sum_{W' \in S} \lambda_{W'} \eCnt{D}{W'}
\]
for some set $S$ of closed walks of $F$ and $\lambda_{W'} \in \mathbb{Q}$ for $W' \in S$, then
\[
\eCnt{D'}{\pi_P(W)} = \sum_{W' \in S} \lambda_{W'} \eCnt{D'}{\pi_P(W')}.
\]
\end{corollary}
\begin{proof}
Let $\mathbf{o}$ be a transitional orientation of $F$ such that $\dTch{P}{\mathbf{o}} = D'$. By Theorem~\ref{thm:pi_on_vectors_full},
\begin{align*}
\eCnt{\dTch{P}{\mathbf{o}}}{\pi_P(W)} 
& = \eCnt{D}{W} \cdot \edgetrans{D}{\mathbf{o}}|_{\tau(P)} \\
& = \left(\sum_{W' \in S} \lambda_{W'} \eCnt{D}{W'}\right) \cdot \edgetrans{D}{\mathbf{o}}|_{\tau(P)} \\
& = \sum_{W' \in S} \lambda_{W'} \eCnt{D}{W'} \cdot \edgetrans{D}{\mathbf{o}}|_{\tau(P)}\\
& = \sum_{W' \in S} \lambda_{W'} \eCnt{\dTch{P}{\mathbf{o}}}{\pi_P(W')},
\end{align*}
where we used Theorem~\ref{thm:pi_on_vectors_full} again in the last equality.
\end{proof}

\begin{lemma}\label{lem:transfer_cycle_spanning}
Let $F$ be a $4$-regular graph and let $P$ be a circuit partition of $F$.

If $\Gamma \supseteq P$ is an (integral, resp.)\ cycle spanning set of $F$, then $\pi_P(\Gamma \setminus P)$ is an (integral, resp.)\ cycle spanning set of $\Tch(P)$.
\end{lemma}
\begin{proof}
This follows from Corollary~\ref{cor:transfer_weights} and the facts that (1) $\pi_P$ maps elements of $P$ to the empty set and (2) every closed walk $W'$ of $\dTch{P}{\mathbf{o}}$ is of the form $W' = \pi_P(W)$ for some closed walk $W$ of $F$.
\end{proof}

We now prove a stronger version of Lemma~\ref{lem:transfer_cycle_spanning}, cf.\ Theorem~\ref{thm:transfer_cycle_spanning_E} below, in order to later prove the strictness condition in Theorem~\ref{thm:main_result}. For a graph $G$ and $E \subseteq E(G)$ we denote by $G-E$ the graph obtained from $G$ by removing the edges of $E$. 

\begin{lemma}\label{lem:cycle_spanning_E_cpartition}
Let $F$ be a $4$-regular graph and let $P$ be a circuit partition of $F$. Let $E \subseteq E(F)$ be such that each circuit of $P$ traverses at most one edge from $E$.

If $\Gamma$ is an (integral, resp.)\ cycle spanning set of $F-E$, then $\Gamma \cup P$ is an (integral, resp.)\ cycle spanning set of $F$.

\end{lemma}
\begin{proof}
Let $\Gamma$ be a cycle spanning set of $F-E$ and let $D$ be a directed version of $F$. Let $W$ be a closed walk of $F$ and let $E' \subseteq E$ be the edges of $E$ that are traversed by $W$. We prove by induction on $|E'|$ that $\eCnt{D}{W} \in \spn_{\mathbb{Q}}(\eCnt{D}{\Gamma \cup P})$.

Assume first that $|E'|=0$, i.e., $E' = \emptyset$. Then $W$ is a closed walk of $F-E$ and so $\eCnt{D}{W} \in \spn_{\mathbb{Q}}(\eCnt{D}{\Gamma})$, since $\Gamma$ is a cycle spanning set of $F-E$. Thus, $\eCnt{D}{W} \in \spn_{\mathbb{Q}}(\eCnt{D}{\Gamma \cup P})$.

Assume now that $|E'|>0$. Let $e \in E'$. Let $p \in P$ be the circuit of $F$ traversing $e$. Consider the closed walk $W'$ obtained from $W$ that avoids the traversal of $e$ by instead taking the path obtained from $p$ by removing $e$. By the induction hypothesis, $\eCnt{D}{W'} \in \spn_{\mathbb{Q}}(\eCnt{D}{\Gamma})$. Now, $\eCnt{D}{W} = \eCnt{D}{W'} + \lambda_{p_d} \cdot \eCnt{D}{p_d}$, where $\lambda_{p_d} \in \{-1,1\}$ and $p_d$ is some orientation of $p$. So, $\eCnt{D}{W} \in \spn_{\mathbb{Q}}(\eCnt{D}{\Gamma \cup P})$.

Since the $\lambda_{p_d}$'s are in $\{-1,1\}$, we have that if $\Gamma$ is integral, then so is $\Gamma \cup P$.
\end{proof}

\begin{theorem}\label{thm:transfer_cycle_spanning_E}
Let $F$ be a $4$-regular graph and let $P$ be a circuit partition of $F$. Let $E \subseteq E(F)$ be such that each circuit of $P$ traverses at most one edge from $E$.

If $\Gamma$ is an (integral, resp.)\ cycle spanning set of $F-E$, then $\pi_P(\Gamma \setminus P)$ is an (integral, resp.)\ cycle spanning set of $\Tch(P)$.
\end{theorem}
\begin{proof}
Let $\Gamma$ be an (integral, resp.)\ cycle spanning set of $F-E$. By Lemma~\ref{lem:cycle_spanning_E_cpartition}, $\Gamma \cup P$ is an (integral, resp.)\ cycle spanning set of $F$. By Lemma~\ref{lem:transfer_cycle_spanning}, $\pi_P((\Gamma \cup P) \setminus P) = \pi_P(\Gamma \setminus P)$ is an (integral, resp.)\ cycle spanning set of $\Tch(P)$.
\end{proof}

\section{Cycle matrices}\label{sec:cycle_matrices}

In this section we prove Theorem~\ref{thm:main_result}. Here the edge-transition incidence matrix introduced in the previous section plays a central role in the construction of the totally transversally unimodular representation of $\mathcal{Z}_3(G)$ for a circle graph $G$.

The \emph{cycle matrix} of a set of closed walks $\Gamma$ of a directed version $D$ of a graph $G$, denoted by $\mathrm{CM}(G,\Gamma,D)$, is the $\Gamma \times E(G)$-matrix over $\mathbb{Z}$ where the row indexed by $W \in \Gamma$ is $\eCnt{D}{W}$. Note that for directed versions $D_1$ and $D_2$ of $G$, the cycle matrix of $\Gamma$ w.r.t.\ $D_1$ is obtained from the cycle matrix of $\Gamma$ w.r.t.\ $D_2$ by multiplying some (possibly none) columns by $-1$. Also note that if $\Gamma$ is a cycle spanning set, then $\mathrm{CM}(G,\Gamma,D)$ is a representation of the cographic matroid $M^*(G)$, that is, the dual of the cycle matroid of $G$. 

By Theorem~\ref{thm:pi_on_vectors_full}, we have the following.
\begin{corollary}\label{cor:cycle_matrix_mult}
Let $F$ be a $4$-regular graph and let $D$ be a directed version of $F$. Let $\mathbf{o}$ be a transitional orientation of $F$ and let $\Gamma$ be a set of closed walks of $D$. Let $P$ be a circuit partition of $F$. Then 
\[
\mathrm{CM}(\Tch(P),\pi_P(\Gamma),\dTch{P}{\mathbf{o}}) \text{ is equal to } \mathrm{CM}(F,\Gamma,D) \cdot \edgetrans{D}{\mathbf{o}}|_{\tau(P)}
\]
up to relabeling of each row index $\pi_P(W)$ of the former matrix by $W$.
\end{corollary}

By Theorem~\ref{thm:transfer_cycle_spanning_E} and Corollary~\ref{cor:cycle_matrix_mult}
we have the following.
\begin{corollary}\label{cor:cs_4reg_touch}
Let $F$ be a $4$-regular graph and let $D$ be a directed version of $F$. Let $\mathbf{o}$ be a transitional orientation of $F$. Let $P$ be a circuit partition of $F$. Let $E \subseteq E(F)$ be such that each circuit of $P$ traverses at most one edge from $E$. Let $\Gamma$ be a cycle spanning set of $F-E$.

Then $\mathrm{CM}(F,\Gamma,D) \cdot \edgetrans{D}{\mathbf{o}}|_{\tau(P)}$ represents $M^*(\Tch(P))$.
\end{corollary}

The \emph{Eulerian} $3$-matroid $Q(F)$ of a $4$-regular graph $F$ is the (unique) $3$-matroid $(\mathfrak{T}(F),\allowbreak\Omega,r)$, where (1) $\Omega = \{ \omega_v \mid v \in V(F) \}$ and, for $v \in V(F)$, $\omega_v$ is the set of transitions at $v$, and (2) for each transversal $T \in \mathcal{T}(\Omega)$, $r(T) = |V(F)| - (|P| - c(F))$ (in other words, the nullity of $T$ is $|P| - c(F)$), where $P$ is the circuit partition with $\tau(P) = T$; see \cite{DBLP:journals/siamdm/Bouchet97}. 

We remark that these conditions uniquely determine $Q(F)$, since the rank function $r$ of any multimatroid where each skew class has cardinality at least two is uniquely determined by the set of transversals $T$ with $r(T)=|T|$; see \cite[Proposition 5.5]{DBLP:journals/siamdm/Bouchet97}.

If $P$ is a circuit partition then the transverse matroid $Q(F)[\tau(P)]$ is equal to $M^*(\Tch(P))$, see \cite[Sec.\ 4]{Bouchet/MMIV} or \cite[Sec.\ 5]{EUJC/Traldi/transmat}.

We thus have the following.
\begin{theorem}\label{thm:repr_3m}
Let $F$ be a $4$-regular graph and let $D$ be a directed version of $F$. Let $\mathbf{o}$ be a transitional orientation of $F$ and let $\Gamma$ be a cycle spanning set of $F - E$, where $E \subseteq E(F)$ contains at most one edge from each connected component of $F$.

Then $\mathrm{CM}(F,\Gamma,D) \cdot \edgetrans{D}{\mathbf{o}}$ represents $Q(F)$. This representation is strict when $|E| = c(F)$.
\end{theorem}
\begin{proof}
Only the last statement is left to show. The rank of $\mathrm{CM}(F,\Gamma,D)$ is equal to the rank of $\mathrm{CM}(F-E,\Gamma,D)$, which in turn is equal to $r(M^*(F-E))$. We have $r(M^*(F-E)) = |E(F-E)| - (|V(F-E)|-c(F-E)) = |E(F)|-|E|-(|V(F)|-c(F)) = 2|V(F)| - |E| - |V(F)| + c(F) = |V(F)| - |E| + c(F)$. Thus, if $|E| = c(F)$, then the rank of $\mathrm{CM}(F,\Gamma,D)$ is $|V(F)|$. Consequently, the rank of the product $\mathrm{CM}(F,\Gamma,D) \cdot \edgetrans{D}{\mathbf{o}}$ is at most $|V(F)|$. 

On the other hand, if $C$ is an Euler system of $F$ then every edge of $\Tch(C)$ is a loop, so the rank of $M^*(\Tch(P))$ is $|V(F)|$. The transverse matroid $Q(F)[\tau(C)]$ is equal to $M^*(\Tch(C))$, so $|V(F)|$ equals the rank of the submatrix of $\mathrm{CM}(F,\Gamma,D) \cdot \edgetrans{D}{\mathbf{o}}$ consisting of columns corresponding to elements of $\tau(C)$. Therefore the rank of $\mathrm{CM}(F,\Gamma,D) \cdot \edgetrans{D}{\mathbf{o}}$ is at least $|V(F)|$. 
\end{proof}

We now recall the following well-known result (see, e.g., \cite[Theorem~19.3]{Schrijver/LIP}).
\begin{proposition}\label{prop:char_tot_unimod}
An $X \times Y$-matrix $A$ is totally unimodular if and only if for every $Z \subseteq X$, there are $\lambda_z \in \{-1,1\}$, for all $z \in Z$, such that all entries of the vector $\sum_{z \in Z} \lambda_z A_{z,\bullet}$ are in $\{-1,0,1\}$, where $A_{z,\bullet}$ is the row vector of $A$ indexed by $z$.
\end{proposition}

Let $G$ be a graph and let $T$ be a maximal forest of $G$. For $e \in E(G) \setminus E(T)$, the unique circuit for which the only edge that is traversed outside $T$ is $e$, is called the \emph{fundamental circuit} of $e$ w.r.t.\ $T$. Let $B$ be the set of oriented circuits obtained by fixing an arbitrary orientation to each fundamental circuit $C_e$ of $e \in E(G) \setminus E(T)$ w.r.t.\ $T$. It is well known that $B$ is an integral cycle basis of $G$ (see, e.g., \cite{CycleBases/DAM/LR}). Let us call $B$ a \emph{strictly fundamental} cycle basis of $G$ w.r.t.\ $T$. The following is well known, see, e.g., \cite[Chapter~19]{Schrijver/LIP}.

\begin{lemma}\label{lem:unique_zero_one_repr}
Let $D$ be a directed version of a graph $G$, let $T$ be a maximal forest of $G$, and let $B$ be a strictly fundamental cycle basis of $G$ w.r.t.\ $T$. Then:
\begin{itemize}
\item $\mathrm{CM}(G,B,D)$ is totally unimodular. 
\item For every oriented cycle $C$ of $G$, $\eCnt{D}{C}$ is equal to $\sum_{C' \in B} \lambda_{C',C} \eCnt{D}{C'}$ with $\lambda_{C',C} \in \{-1,0,1\}$ for all $C' \in B$.
\item If $C$ is an oriented cycle of $G$ and $C_e \in B$ denotes an oriented fundamental circuit for $e \in E(G) \setminus E(T)$ with respect to $T$, then $\lambda_{C_e,C} \neq 0$ if and only if $e$ is traversed by $C$.
\end{itemize}
\end{lemma}

The second property of Lemma~\ref{lem:unique_zero_one_repr} is somewhat similar to the notion of a zero-one cycle basis considered in \cite{CycleBasesSurvey/KavithaLMMRUZ09}. However, here we allow vertex repetitions in oriented circuits, cf.\ Remark~\ref{rem:bases_vertex_rep}. Example~\ref{ex:TU_not_implies_TU} below illustrates that this is a crucial difference.

\begin{lemma}\label{lem:circuit_TchP}
Let $F$ be a $4$-regular graph and let $P$ be a circuit partition of $F$. For each cycle $C$ of $F$, there is a cycle $C'$ of $F$ that traverses the same set of edges as $C$ such that $\pi_P(C')$ is a cycle of $\Tch(P)$.
\end{lemma}
\begin{proof}
Let $C'$ be a cycle of $F$ obtained from $C$ by changing, for each vertex $v$ of $F$ for which all four incident half-edges are traversed by some circuits of $C$, the transition $t$ taken by $C$ at $v$ such that it coincides with the transition taken by $P$ at $v$ (of course, $t$ might already coincide with the transition of $P$ at $v$, in which case we change nothing on $C$ at $v$). By the construction, $\pi_P(C')$ does not traverse any edge of $\Tch(P)$ twice, so $\pi_P(C')$ is a cycle of $\Tch(P)$.
\end{proof}

\begin{theorem}\label{thm:Tch_totunimod}
Let $F$ be a $4$-regular graph, let $P$ be a circuit partition of $F$, and let $D$ be a directed version of $\Tch(P)$. Let $B$ be a strictly fundamental cycle basis of $F$.

Then $\mathrm{CM}(\Tch(P),\pi_P(B),D)$ is totally unimodular.
\end{theorem}
\begin{proof}
Let $Z \subseteq \pi_P(B)$. By Proposition~\ref{prop:char_tot_unimod}, it suffices to show that there are $\lambda_z \in \{-1,1\}$, for all $z \in Z$, such that all entries of the vector $\sum_{z \in Z} \lambda_z \eCnt{D}{z}$ are in $\{-1,0,1\}$. Let $Z' \subseteq B$ such that $\pi_P(Z') = Z$. Since $\mathrm{CM}(F,B,D)$ is totally unimodular, by Proposition~\ref{prop:char_tot_unimod} we (may) choose some assignment of $\lambda'_z \in \{-1,1\}$ for all $z \in Z'$ such that all entries of $\sum_{z \in Z'} \lambda'_{z} \eCnt{D}{z}$ are in $\{-1,0,1\}$. By Lemma~\ref{lem:sum_oriented_cycles}, there is an oriented cycle $C$ of $F$ such that $\eCnt{D}{C} = \sum_{z \in Z'} \lambda'_{z} \eCnt{D}{z}$. By Lemma~\ref{lem:circuit_TchP}, there is an oriented cycle $C'$ of $F$ such that $C$ and $C'$ traverse the same edges and $\pi_P(C')$ is an oriented cycle of $\Tch(P)$. By Lemma~\ref{lem:unique_zero_one_repr}, $\eCnt{D}{C'} = \sum_{z \in Z'} \lambda''_z \eCnt{D}{z}$ with $\lambda''_z \in \{-1,1\}$ for all $z \in Z'$. By Corollary~\ref{cor:transfer_weights}, $\eCnt{D}{\pi_P(C')}$ is equal to $\sum_{z \in Z'} \lambda''_z \eCnt{D}{\pi_P(z)}$, which in turn is equal to $\sum_{z \in Z} \lambda_{z} \eCnt{D}{z}$ by setting $\lambda_{\pi_P(z)} := \lambda''_z$ for all $z \in Z'$. Thus all entries of $\sum_{z \in Z} \lambda_z \eCnt{D}{z}$ are indeed in $\{-1,0,1\}$.
\end{proof}

For an $X \times \mathfrak{T}(F)$-matrix $A$, where $F$ is a $4$-regular graph, we say that $A$ is \emph{totally transversally unimodular} if for every transversal $T$ of $\mathfrak{T}(F)$, the $X \times T$-submatrix of $A$ induced by the columns of $T$ is totally unimodular.

By Corollary~\ref{cor:cycle_matrix_mult} and Theorem~\ref{thm:Tch_totunimod} we obtain the following.
\begin{corollary}\label{cor:trans_totunimod}
Let $F$ be a $4$-regular graph, let $D$ be a directed version of $F$, let $B$ be a strictly fundamental cycle basis of $F$, and let $\mathbf{o}$ be a transitional orientation of $F$.

Then $\mathrm{CM}(F,B,D) \cdot \edgetrans{D}{\mathbf{o}}$ is totally transversally unimodular.
\end{corollary}

We now obtain the following.
\begin{theorem}\label{thm:eulerian_3m_ttu}
Let $F$ be a $4$-regular graph. Then $Q(F)$ has a strict, totally transversally unimodular representation.
\end{theorem}
\begin{proof}
Let $B$ be a strictly fundamental cycle basis of $F - E$, where $E$ contains exactly one edge of each connected component of $F$. By Theorem~\ref{thm:repr_3m}, $A := \mathrm{CM}(F,B,D) \cdot \edgetrans{D}{\mathbf{o}}$ is a strict representation of $Q(G)$. Let $B'$ be the strictly fundamental cycle basis of $F$ such that $B \subseteq B'$. By Corollary~\ref{cor:trans_totunimod}, $A' := \mathrm{CM}(F,B',D) \cdot \edgetrans{D}{\mathbf{o}}$ is totally transversally unimodular. Since $A$ is obtained from $A'$ by removing some rows, $A$ is totally transversally unimodular too.
\end{proof}

An \emph{interlacement graph} $G$ of a $4$-regular graph $F$ with respect to some Euler system $C$ is a simple graph without loops such that $V(G) = V(F)$ and there is an edge between distinct vertices $u$ and $v$ if and only if $u$ and $v$ belong to the same connected component of $F$ and the Eulerian circuit $C'$ of $C$ corresponding to that connected component visits $u$ and $v$ in the order $u,v,u,v$ or $v,u,v,u$ (i.e., $u$ and $v$ are ``interlaced'' in $C$). We recall that a \emph{circle graph} is an interlacement graph of some $4$-regular graph with respect to some Euler system. 

\begin{proposition}[\cite{EUJC/Traldi/transmat}]\label{prop:circle_graph_3m_Eulerian_3m}
Let $G$ be the interlacement graph of a $4$-regular graph $F$ with respect to some Euler system $C$. Then $\mathcal{Z}_3(G)$ is equal to $Q(F)$ up to relabelling, for each vertex $v \in V(F)$, $\phi_C(v)$, $\chi_C(v)$, and $\psi_C(v)$ to the column indices corresponding to the $v$-columns of $I$, $A(G)$, and $I+A(G)$ in $\mathrm{IAS}(G)$, respectively.
\end{proposition}

By Theorem~\ref{thm:eulerian_3m_ttu} and Proposition~\ref{prop:circle_graph_3m_Eulerian_3m} we obtain the main result of this paper (cf.\ Theorem~\ref{thm:main_result} in the introduction).
\begin{theorem}\label{thm:circle_graph_ttu}
Let $G$ be a circle graph. Then $\mathcal{Z}_3(G)$ has a strict, totally transversally unimodular representation.
\end{theorem}

\begin{remark} \label{remark:no_TU_repr}
In view of Theorem~\ref{thm:circle_graph_ttu} one may wonder whether $\mathcal{Z}_3(G)$ has a totally  unimodular representation (i.e., is sheltered by some regular matroid) for circle graphs $G$. This turns out to be false in general. Indeed, consider the circle graph $G = K_4$. It is easy to see that the matroid represented by $(I \enskip A(K_4))$ is not regular (indeed, the representation of the dual Fano matroid $F_7^*$ of the introduction is a submatrix of this matrix). By \cite[Proposition 18]{BT/IsotropicMatI}, for any simple graph $G$, the 2-matroid $\mathcal{Z}_3(G) - \Psi$, where $\Psi$ is the set of column indices corresponding to $I+A(G)$ in $\mathrm{IAS}(G)$, is ``tight'' (meaning that the union of every full rank transversal and every skew class has a circuit of a transverse matroid as a subset). By \cite[Propositions 17 and 24]{BT/IsotropicMatI}, for any tight 2-matroid $Z$ there is at most one binary matroid $M$ that shelters $Z$. Therefore, no regular matroid shelters $\mathcal{Z}_3(K_4) - \Psi$, and so no regular matroid shelters $\mathcal{Z}_3(K_4)$.
\end{remark}

\begin{example}
Consider again $F$, $C$, and $D$ from the running example, see Figure~\ref{fig:4reg_graph} (for $F$) and Figure~\ref{fig:Euler} (for $C$ and $D$). Let $T$ be the spanning tree of $F$ consisting of edges $e_1$, $e_2$, and $e_3$. Let $B$ be a strictly fundamental cycle basis of $F$ with respect to $T$, where each oriented fundamental circuit $C_e$ of $e$ is oriented such that $e$ is traversed in the direction of $D$. The totally unimodular matrix $\mathrm{CM}(F,B,D)$ is as follows:
\[
\bordermatrix{
~   & e_1 & e_2 & e_3 & e_4 & e_5 & e_6 & e_7 & e_8 \cr
C_{e_4} &  0 &  1 &  1 &  1 &  0 &  0 &  0 &  0 \cr
C_{e_5} &  1 &  0 &  0 &  0 &  1 &  0 &  0 &  0 \cr
C_{e_6} & -1 & -1 &  0 &  0 &  0 &  1 &  0 &  0 \cr
C_{e_7} &  0 &  0 & -1 &  0 &  0 &  0 &  1 &  0 \cr
C_{e_8} &  1 &  1 &  1 &  0 &  0 &  0 &  0 &  1
}.
\]

The matrix $\mathrm{CM}(F,B,D) \cdot \edgetrans{D}{\mathbf{o}}$ is equal to
\[
\scalebox{0.9}{
{\let\quad\thinspace
\bordermatrix{
~ & 
\phi_C(a) & \phi_C(b) & \phi_C(c) & \phi_C(d) &  
\chi_C(a) & \chi_C(b) & \chi_C(c) & \chi_C(d) & 
\psi_C(a) & \psi_C(b) & \psi_C(c) & \psi_C(d) \cr
C_{e_4} &  0 &  1 &  0 &  0 &  0 &  0 & -1 & -1 &  0 &  1 &  1 &  1 \cr
C_{e_5} &  1 & -1 &  0 &  0 &  0 &  0 &  0 &  0 &  1 &  1 &  0 &  0 \cr
C_{e_6} & -1 &  0 &  1 &  0 &  1 &  1 &  1 &  0 &  0 & -1 &  0 &  0 \cr
C_{e_7} &  0 &  0 & -1 &  1 &  0 &  0 &  1 &  1 &  0 &  0 &  0 &  0 \cr
C_{e_8} &  0 &  0 &  0 & -1 & -1 & -1 & -1 &  0 &  1 &  1 &  1 &  1
}}}.
\]
By Theorem~\ref{thm:repr_3m} (with $E = \emptyset$), this matrix represents $Q(F)$. In fact, taking $E$ to be a singleton instead of the empty set, we have by Theorem~\ref{thm:repr_3m} that the matrix obtained by removing any row of this matrix remains a representation of $Q(F)$. By Corollary~\ref{cor:trans_totunimod}, the above matrix is totally transversally unimodular.
\end{example}

One may wonder, in view of Corollary~\ref{cor:trans_totunimod}, whether total unimodularity of $\mathrm{CM}(F,B,D)$ directly implies total transversal unimodularity of $\mathrm{CM}(F,B,D) \cdot \edgetrans{D}{\mathbf{o}}$. The next example shows that this is not the case.

\begin{example}\label{ex:TU_not_implies_TU}
Let $F$ be a $4$-regular graph with two vertices and four parallel edges between the vertices and let $D$ be a directed version of $F$ such that each vertex has two incoming and two outgoing edges. Let $(e_1, e_2, e_3, e_4)$ be a closed walk representing an Eulerian circuit $C_1$ of $D$ (for notational convenience we represent closed walks in this example as sequences of edges instead of directed single transitions). Let $B$ be the cycle basis of $F$ consisting of the oriented circuits $C_{1,B}$, $C_{2,B}$, and $C_{3,B}$ corresponding to the closed walks $(e_1,e_2)$, $(e_2,e_3)$, and $(e_3,e_4)$, respectively. Consider the oriented circuit $C_2$ corresponding to the closed walk $(e_1,e_4,\bar e_2,\bar e_3)$, where $\bar e$ again means traversing the directed edge $e$ in the opposite direction. While $C_2$ is an oriented cycle, one verifies that it cannot be written as $\sum_{C \in B} \lambda_{C} \eCnt{D}{C}$ where $\lambda_{C} \in \{-1,0,1\}$ for all $C \in B$ (note that it can be written in such a way if $B$ is replaced by a strictly fundamental cycle basis, cf.\ Lemma~\ref{lem:unique_zero_one_repr}). However, $B$ is an integral cycle basis since (1) for every oriented circuit $C$, $\eCnt{D}{C}$ is the sum of $\eCnt{D}{C_i}$'s, where each $C_i$ is an oriented circuit without vertex repetitions, and (2) one easily verifies that, for every oriented circuit $C_i$ without vertex repetitions, $\eCnt{D}{C_i}$ can be written as the integral sum of the $\eCnt{D}{C'}$'s where $C' \in B$.

It is interesting to observe that every oriented circuit that visits each vertex at most once \emph{can} be written as $\sum_{C \in B} \lambda_{C} \eCnt{D}{C}$ where $\lambda_{C} \in \{-1,0,1\}$ for all $C \in B$. This property of $B$ is captured by the notion of a zero-one cycle basis considered in \cite{CycleBasesSurvey/KavithaLMMRUZ09}. Therefore, this property is crucially different from the second property given in Lemma~\ref{lem:unique_zero_one_repr}.

Let $\mathbf{o}$ be the transitional orientation of $F$ such that $\mathbf{o}(t)$ always chooses the single transition containing the half-edge among the half-edges incident to the vertex corresponding to $t$ that is traversed first by $C_1$ starting half-way in $e_1$ and walking in the direction of $e_1$. We have that
\[
\mathrm{CM}(F,B,D) = 
\bordermatrix{
~  & e_1 & e_2 & e_3 & e_4\cr
C_{1,B} &  1 &  1 &  0 &  0\cr
C_{2,B} &  0 &  1 &  1 &  0\cr
C_{3,B} &  0 &  0 &  1 &  1
}
\]
is totally unimodular and
\[
\edgetrans{D}{\mathbf{o}} = 
\bordermatrix{
~  & \phi_{C_1}(v_1) & \phi_{C_1}(v_2) & \chi_{C_1}(v_1) & \chi_{C_1}(v_2) & \psi_{C_1}(v_1) & \psi_{C_1}(v_2)\cr
e_1 & -1 &  0 & -1 &  1 & -1 &  0\cr
e_2 &  1 & -1 &  0 & -1 &  0 & -1\cr
e_3 &  0 &  1 &  0 &  0 & -1 &  0\cr
e_4 &  0 &  0 &  1 &  0 &  0 & -1
}.
\]
The matrix $A = \mathrm{CM}(F,B,D) \cdot \edgetrans{D}{\mathbf{o}}$ is equal to
\[
\bordermatrix{
~  & \phi_{C_1}(v_1) & \phi_{C_1}(v_2) & \chi_{C_1}(v_1) & \chi_{C_1}(v_2) & \psi_{C_1}(v_1) & \psi_{C_1}(v_2)\cr
C_{1,B} &  0 & -1 & -1 &  0 & -1 & -1\cr
C_{2,B} &  1 &  0 &  0 & -1 & -1 & -1\cr
C_{3,B} &  0 &  1 &  1 &  0 & -1 & -1\cr
}.
\]
By Corollary~\ref{cor:cycle_matrix_mult}, $\mathrm{CM}(\Tch(P),\pi_P(B),\dTch{P}{\mathbf{o}}) = A|_{\tau(P)}$ for every circuit partition $P$ of $F$. Consider the circuit partition $P$ with $\tau(P) = \{\chi_{C_1}(v_1),\psi_{C_1}(v_2)\}$. Then $A|_{\tau(P)}$ is not totally unimodular. Indeed, its submatrix induced by the rows indexed by $C_{1,B}$ and $C_{3,B}$ has determinant $2$.

Alternatively, if we take the strictly fundamental cycle basis $B'$ of $G$ w.r.t.\ the spanning tree $T$ consisting of edge $e_1$ such that each oriented circuit in $B'$ is oriented in the direction of $e_1$, then we get
\[
\mathrm{CM}(F,B',D) = 
\bordermatrix{
~  & e_1 & e_2 & e_3 & e_4\cr
C_{1,B'} &  1 &  1 &  0 &  0\cr
C_{2,B'} &  1 &  0 & -1 &  0\cr
C_{3,B'} &  1 &  0 &  0 &  1
}
\]
and $A' = \mathrm{CM}(F,B',D) \cdot \edgetrans{D}{\mathbf{o}}$ is equal to
\[
\bordermatrix{
~  & \phi_{C_1}(v_1) & \phi_{C_1}(v_2) & \chi_{C_1}(v_1) & \chi_{C_1}(v_2) & \psi_{C_1}(v_1) & \psi_{C_1}(v_2)\cr
C_{1,B'} &  0 & -1 & -1 &  0 & -1 & -1\cr
C_{2,B'} & -1 & -1 & -1 &  1 &  0 &  0\cr
C_{3,B'} & -1 &  0 &  0 &  1 & -1 & -1\cr
}.
\]
By Corollary~\ref{cor:trans_totunimod}, $A'$ is transversally totally unimodular.
\end{example}

\section{Circuits induced by an Eulerian circuit}\label{sec:cycle_basis_Euler}

In this section we discuss how the totally transversally unimodular representations of $\mathcal{Z}_3(G)$ considered in this paper relate to the representations considered in \cite{BT/IsotropicMatIV}. In this way we obtain simpler proofs of some of the results of \cite{BT/IsotropicMatIV}.

First we recall the following result from \cite{Liebchen03/CycleBases}. For convenience we also give a proof. Recall that a square matrix of integers is \emph{unimodular} if its determinant is $\pm 1$.
\begin{proposition}[\cite{Liebchen03/CycleBases}]\label{prop:cycle_spanning_unimodular}
Let $G$ be a graph, $B$ be an integral cycle basis of $G$, and $D$ be some directed version of $G$. Let $T$ be a set of edges of $G$ that forms a maximal forest of $G$. Then the matrix obtained from $\mathrm{CM}(G,B,D)$ by removing the columns of $T$ is unimodular.
\end{proposition}
\begin{proof}
Since $B$ is an integral cycle basis, we have $\eCnt{D}{C} \in \mathrm{span}_{\mathbb{Z}}(\eCnt{D}{B})$ for every oriented circuit $C$. Let $e \in E(D)\setminus T$. For the oriented fundamental circuit $C_e$ of $e$ traversing in the direction of $e$ in $D$, the restriction of $\eCnt{D}{C_e}$ to index set $E(D)\setminus T$ is a unit vector with the entry of $e$ equal to $1$. Consequently, the span over $\mathbb{Z}$ of the rows of the matrix $A$ obtained from $\mathrm{CM}(G,B,D)$ by removing the columns of $T$ is $\mathbb{Z}^{E(D)\setminus T}$. By \cite[Theorem~4.3]{Schrijver/LIP}, $A$ is unimodular.
\end{proof}

For a graph $G$, we say that $E \subseteq E(G)$ is \emph{based} in $G$ if $E$ contains exactly one edge of each connected component of $G$. We now provide a counterpart to Theorem~\ref{thm:Tch_totunimod}.
\begin{theorem}\label{thm:Tch_unimod}
Let $F$ be a $4$-regular graph, let $P$ be a circuit partition of $F$, and let $D$ be a directed version of $\Tch(P)$. Let $\Gamma$ be an integral cycle basis of $F-E$, where $E \subseteq E(F)$ is based in $F$.

Then $\mathrm{CM}(\Tch(P),\pi_P(\Gamma),D)$ is square and has determinant $-1$, $0$, or $1$.
\end{theorem}
\begin{proof}
Let $\Gamma$ be an integral cycle basis of $F-E$. Then $\mathrm{CM}(\Tch(P),\pi_P(\Gamma),D)$ has $|E(\Tch(P))| = |V(F)|$ columns and $|\pi_P(\Gamma)| = |\Gamma|$ rows. Since $\Gamma$ is a cycle basis of $F-E$, we have $|\Gamma| = |E(F-E)| - (|V(F-E)| - c(F-E)) = |E(F)| - c(F) - (|V(F)| - c(F)) = |E(F)| - |V(F)| = |V(F)|$, and so the matrix is indeed square. 

If $C \in \Gamma \cap P$, then the row of $\mathrm{CM}(\Tch(P),\pi_P(\Gamma),D)$ indexed by $C$ is zero, so the statement holds. Assume now that $\Gamma \cap P = \emptyset$. By Theorem~\ref{thm:transfer_cycle_spanning_E}, $\pi_P(\Gamma)$ is an integral cycle spanning set of $\Tch(P)$. If $\det(\mathrm{CM}(\Tch(P),\pi_P(\Gamma),D)) \neq 0$, then $\pi_P(\Gamma)$ is an integral cycle basis of $\Tch(P)$. By Proposition~\ref{prop:cycle_spanning_unimodular}, we obtain that $\mathrm{CM}(\Tch(P),\pi_P(\Gamma),D)$ is unimodular (and every edge of $\Tch(P)$ is a loop).
\end{proof}

The following (easy to verify) result is used in the BEST theorem \cite{BEST_thm_first_half,BEST_thm_second_half}.
\begin{lemma}[\cite{BEST_thm_first_half,BEST_thm_second_half}]\label{lem:BEST_spanning_tree}
Let $F$ be a connected $4$-regular graph and let $C$ be an oriented Eulerian circuit of $F$. Let $e$ be an edge of $F$. Let $T$ be the graph obtained from $F$ by removing $e$ and removing, for every $v \in V(F)$, the incoming edge $R_{C,e}(v)$ when visiting $v$ for the second time while walking along $C$ and starting at the middle of $e$. Then $T$ is a spanning tree of $F$.
\end{lemma}

We denote the half-edge of $R_{C,e}(v)$ (from Lemma~\ref{lem:BEST_spanning_tree}) incident to $v$ by $H_{C,e}(v)$.

The spanning tree of Lemma~\ref{lem:BEST_spanning_tree} is called the spanning tree of $F$ \emph{induced} by $C$ and $e$.

Obviously, we can apply the above lemma to each connected component of a $4$-regular graph. So, if $F$ is a $4$-regular graph, $C$ is an oriented Euler system of $F$, and $E \subseteq E(F)$ is based in $F$, then we (may) speak of the \emph{maximal forest} of $F$ induced by $C$ and $E$. Similarly, we define $R_{C,E}(v)$ and $H_{C,E}(v)$ in this more general context.

\begin{definition}
Let $C$ be an oriented Euler system of a $4$-regular graph $F$ and let $E \subseteq E(F)$ be based in $F$. For a vertex $v$ of $F$, the oriented circuit \emph{induced by} $C$ at $v$ based on $E$ is the oriented circuit that traverses the segment from $v$ to $v$ of an oriented Euler circuit of $C$, and avoids traversing edges of $E$.
\end{definition}

We denote by $\Gamma_{E,C}$ the set of all oriented circuits induced by $C$ and based on $E$. Note that the orientations of the oriented circuits of $\Gamma_{E,C}$ coincide with the oriented circuits of $C$. Consequently, if $D$ is a directed version of $F$ and $C,C' \in \Gamma_{E,C}$, and $e$ is an index for which both its entry in $\eCnt{D}{C}$ is nonzero and its entry in $\eCnt{D}{C'}$ is nonzero, then these entries are equal.

\begin{lemma}\label{lem:fund_based_weakly_fund}
Let $F$ be a $4$-regular graph. Let $E \subseteq E(F)$ be based in $F$ and let $C$ be an oriented Euler system of $F$.

Then $\Gamma_{E,C}$ is an integral cycle basis of $F-E$.
\end{lemma}
\begin{proof}
Let $T$ be the maximal forest of $F$ induced by $C$ and $E$. Let $D$ be a directed version of $F$. For $v \in V(F)$, let $C_v \in \Gamma_{E,C}$ be the oriented circuit induced by $C$ at $v$ based on $E$.

Let $v_1,\ldots,v_n$ be a linear ordering of the vertices of $F$ such that if edge $R_{C,E}(v_i)$ is traversed before $R_{C,E}(v_j)$ in some Eulerian circuit of $C$ starting from some $e \in E$ in the direction coinciding with $D$, then $i < j$.

Notice that for all $i \in \{1,\ldots,n\}$, edge $R_{C,E}(v_i)$ is traversed by $C_{v_i}$, but not by any $C_{v_k}$, $k \in \{1,\ldots,i-1\}$. Moreover, each edge traversed by $C_{v_i}$ outside $T$ is of the form $R_{C,E}(v_k)$ for some $k \in \{1,\ldots,i-1\}$. Hence by substracting appropriate $\eCnt{D}{C_{v_k}}$'s with $k \in \{1,\ldots,i-1\}$ from $\eCnt{D}{C_{v_i}}$ we obtain an element $s$ of the cycle space of $D$ for which every entry is in $\{-1,0,1\}$. By Lemma~\ref{lem:sum_oriented_cycles}, $s$ is the incidence vector of some oriented cycle. Since the only nonzero entry of $s$ outside $T$ is indexed by $R_{C,E}(v_i)$, we observe that $s$ is the incidence vector of the oriented fundamental circuit for $R_{C,E}(v_i)$ with respect to $T$, oriented in the direction of $E$ in $C$. Since the oriented fundamental circuits for $R_{C,E}(v)$ form an integral cycle basis of $F-E$ (because every edge of $F$ outside $T$ and $E$ is of the form $R_{C,E}(v)$), so do the elements $C_v \in \Gamma_{E,C}$.
\end{proof}

\begin{remark}
A cycle basis $B$ of $G$ is called \emph{weakly fundamental}, see, e.g., \cite{CycleBases/DAM/LR}, if there is a linear ordering $(C_1,\ldots,C_n)$ of $B$ such that, for all $i \in \{1,\ldots,n\}$, $C_i$ traverses an edge that is not traversed by any $C_k$ with $k \in \{1,\ldots,i-1\}$. Obviously, every strictly fundamental cycle basis is weakly fundamental. From the proof of Lemma~\ref{lem:fund_based_weakly_fund} we see that the cycle basis $\Gamma_{E,C}$ is weakly fundamental.
\end{remark}

By Lemma~\ref{lem:transfer_cycle_spanning} and Lemma~\ref{lem:fund_based_weakly_fund} we have the following.
\begin{theorem}\label{thm:cycle_spanning_Tch}
Let $F$ be a $4$-regular graph. Let $E \subseteq E(F)$ be based in $F$. Let $P$ be a circuit partition of $F$. Let $C$ be an Euler system of $F$.

Then $\pi_P(\Gamma_{E,C})$ is an integral cycle spanning set of $\Tch(P)$.
\end{theorem}

The next corollary is shown in \cite{BT/IsotropicMatIV} (see Remark~\ref{rem:literature_BT}). In this paper it follows from Theorems~\ref{thm:Tch_unimod} and \ref{thm:cycle_spanning_Tch}.
\begin{corollary}[\cite{BT/IsotropicMatIV}]\label{cor:Tch_unimod}
Let $F$ be a $4$-regular graph. Let $E \subseteq E(F)$ be based in $F$. Let $P$ be a circuit partition of $F$. Let $C$ be an Euler system of $F$. Let $D$ be a directed version of $\Tch(P)$.

Then $\mathrm{CM}(\Tch(P),\pi_P(\Gamma_{E,C}),D)$ has determinant $-1$, $0$, or $1$.
\end{corollary}

Just like Corollary~\ref{cor:trans_totunimod}, Corollary~\ref{cor:Tch_unimod} can be stated independently of the circuit partition $P$ as follows.

\begin{corollary}\label{cor:Tch_unimod_full}
Let $F$ be a $4$-regular graph. Let $E \subseteq E(F)$ be based in $F$. Let $C$ be an Euler system of $F$. Let $D$ be a directed version of $F$ and let $\mathbf{o}$ be a transitional orientation of $F$.

For each transversal $T$ of $\mathfrak{T}(F)$, the submatrix of $\mathrm{CM}(F,\Gamma_{E,C},D) \cdot \edgetrans{D}{\mathbf{o}}$ induced by the columns of $T$ has determinant $-1$, $0$, or $1$.
\end{corollary}

Note that by Corollary~\ref{cor:cycle_matrix_mult}, the matrix $\mathrm{CM}(F,\Gamma_{E,C},D) \cdot \edgetrans{D}{\mathbf{o}}$ given in Corollary~\ref{cor:Tch_unimod_full} does not depend on $D$. If $\mathbf{o}$ is the transitional orientation that assigns to a transition at $v \in V(F)$ the single transition that does not contain the $H_{C,E}(v)$ half-edge, then we denote this matrix by $\mathrm{IAS}(F,C,E)$.

\begin{remark}\label{rem:literature_BT}
In \cite{BT/IsotropicMatIV}, the matrix $M_{\mathbb{R},\Gamma_{E,C}}(C,P,D) := \mathrm{CM}(\Tch(P),\pi_P(\Gamma_{E,C}),D)$ is considered 
and it is shown there (1) that, for any field $\mathbb{F}$, the $\mathbb{F}$-cycle space of $D$ is equal to the $\mathbb{F}$-span of the rows of $M_{\mathbb{R},\Gamma_{E,C}}(C,P,D)$ \cite[Theorem~34]{BT/IsotropicMatIV}, where the $\mathbb{F}$-cycle space and $\mathbb{F}$-span is the ``$\mathbb{F}$-counterpart'' of the cycle space (i.e., over $\mathbb{Q}$) and the integral span (i.e., over $\mathbb{Z}$), respectively, and (2) that $\det(M_{\mathbb{R},\Gamma_{E,C}}(C,P,D)) \allowbreak \in \{-1,0,1\}$ \cite[Corollary~33]{BT/IsotropicMatIV}. Notice that (1) follows from Theorem~\ref{thm:cycle_spanning_Tch} and that (2) follows from Corollary~\ref{cor:Tch_unimod}. These results are shown in \cite{BT/IsotropicMatIV} using a result whose proof relies on examining a large number of different cases separately. The theory developed above allows for alternative and shorter proofs that more deeply explain why these results hold. By the above, one can observe that $\mathrm{IAS}(F,C,E)$ is equal to the matrix $\mathrm{IAS}_{\Gamma^{o}_{E}}(C)$ defined in \cite{BT/IsotropicMatIV}.
\end{remark}

\begin{example}\label{ex:IAS_from_Euler}
Consider again $F$, $D$, $C$, and $\mathbf{o}$ from the running example. The matrix $\mathrm{CM}(F,\Gamma_{E,C},D)$ with $E = \{e_8\}$ is as follows:
\[
\bordermatrix{
~   & e_1 & e_2 & e_3 & e_4 & e_5 & e_6 & e_7 & e_8 \cr
C_a &  1 &  1 &  1 &  1 &  1 &  0 &  0 &  0 \cr
C_b &  0 &  1 &  1 &  1 &  0 &  0 &  0 &  0 \cr
C_c &  0 &  0 &  1 &  1 &  1 &  1 &  0 &  0 \cr
C_d &  0 &  0 &  0 &  1 &  1 &  1 &  1 &  0
},
\]
where $C_v$ is the oriented circuit induced by $C$ at $v$ based on $E$.

Since the columns in the above depiction of $\mathrm{CM}(F,\Gamma_{E,C},D)$ are given in the order of the edges that are visited by $C$ starting at the middle of $e_8$, the row indexed by a $C_v$ consists of the block of $1$'s starting from the first column indexed by an edge having $v$ as its tail until, and including, the second column indexed by an edge $e$ having $v$ as its head (which is edge $R_{C,e}(v)$). The entries outside this block are zero.

Note that $\mathbf{o}$ as defined assigns to a transition at $v \in V(F)$ the single transition that does not contain the $H_{C,E}(v)$ half-edge. Therefore $\mathrm{CM}(F,\Gamma_{E,C},D) \cdot \edgetrans{D}{\mathbf{o}}$ equals $\mathrm{IAS}(F, C, E)$, which in turn is equal to
\[
\scalebox{0.9}{
{\let\quad\thinspace
\bordermatrix{
~ & 
\phi_C(a) & \phi_C(b) & \phi_C(c) & \phi_C(d) &  
\chi_C(a) & \chi_C(b) & \chi_C(c) & \chi_C(d) & 
\psi_C(a) & \psi_C(b) & \psi_C(c) & \psi_C(d) \cr
C_a &  1 &  0 &  0 &  0 &  0 &  0 & -1 & -1 &  1 &  2 &  1 &  1 \cr
C_b &  0 &  1 &  0 &  0 &  0 &  0 & -1 & -1 &  0 &  1 &  1 &  1 \cr
C_c &  0 &  0 &  1 &  0 &  1 &  1 &  0 & -1 &  1 &  1 &  1 &  1 \cr
C_d &  0 &  0 &  0 &  1 &  1 &  1 &  1 &  0 &  1 &  1 &  1 &  1
}}}.
\]
\end{example}

\subsection*{Acknowledgements}
We thank the anonymous referees for helpful comments and corrections, which significantly improved the paper. This research was performed while R.B.\ was a postdoctoral fellow of the Research Foundation -- Flanders (FWO).

\begingroup
\setlength{\emergencystretch}{8em}
\printbibliography
\endgroup

\end{document}